\documentclass[a4paper,10pt]{amsart}
\usepackage[utf8]{inputenc}
\usepackage{amssymb}
\usepackage{amscd}
\usepackage{amsfonts}
\usepackage{amsmath}
\usepackage{amsthm}
\usepackage{mathscinet}
\usepackage{enumerate}
\usepackage{calc}
\usepackage{color}
\usepackage[usenames,dvipsnames,svgnames,table]{xcolor}
\usepackage{environ}
\usepackage{verbatim}
\usepackage{graphicx}
\usepackage{ifthen}
\usepackage{fancyvrb}

\makeatletter{}\makeatletter
\def\PY@reset{\let\PY@it=\relax \let\PY@bf=\relax%
    \let\PY@ul=\relax \let\PY@tc=\relax%
    \let\PY@bc=\relax \let\PY@ff=\relax}
\def\PY@tok#1{\csname PY@tok@#1\endcsname}
\def\PY@toks#1+{\ifx\relax#1\empty\else%
    \PY@tok{#1}\expandafter\PY@toks\fi}
\def\PY@do#1{\PY@bc{\PY@tc{\PY@ul{%
    \PY@it{\PY@bf{\PY@ff{#1}}}}}}}
\def\PY#1#2{\PY@reset\PY@toks#1+\relax+\PY@do{#2}}

\expandafter\def\csname PY@tok@kr\endcsname{\let\PY@bf=\textbf}
\expandafter\def\csname PY@tok@sx\endcsname{\let\PY@it=\textit}
\expandafter\def\csname PY@tok@ss\endcsname{\let\PY@it=\textit}
\expandafter\def\csname PY@tok@sr\endcsname{\let\PY@it=\textit}
\expandafter\def\csname PY@tok@kn\endcsname{\let\PY@bf=\textbf}
\expandafter\def\csname PY@tok@gh\endcsname{\let\PY@bf=\textbf}
\expandafter\def\csname PY@tok@kc\endcsname{\let\PY@bf=\textbf}
\expandafter\def\csname PY@tok@sb\endcsname{\let\PY@it=\textit}
\expandafter\def\csname PY@tok@cs\endcsname{\let\PY@it=\textit}
\expandafter\def\csname PY@tok@ge\endcsname{\let\PY@it=\textit}
\expandafter\def\csname PY@tok@sd\endcsname{\let\PY@it=\textit}
\expandafter\def\csname PY@tok@k\endcsname{\let\PY@bf=\textbf}
\expandafter\def\csname PY@tok@nc\endcsname{\let\PY@bf=\textbf}
\expandafter\def\csname PY@tok@nn\endcsname{\let\PY@bf=\textbf}
\expandafter\def\csname PY@tok@ne\endcsname{\let\PY@bf=\textbf}
\expandafter\def\csname PY@tok@gu\endcsname{\let\PY@bf=\textbf}
\expandafter\def\csname PY@tok@s2\endcsname{\let\PY@it=\textit}
\expandafter\def\csname PY@tok@c1\endcsname{\let\PY@it=\textit}
\expandafter\def\csname PY@tok@gp\endcsname{\let\PY@bf=\textbf}
\expandafter\def\csname PY@tok@c\endcsname{\let\PY@it=\textit}
\expandafter\def\csname PY@tok@cm\endcsname{\let\PY@it=\textit}
\expandafter\def\csname PY@tok@s\endcsname{\let\PY@it=\textit}
\expandafter\def\csname PY@tok@gs\endcsname{\let\PY@bf=\textbf}
\expandafter\def\csname PY@tok@se\endcsname{\let\PY@bf=\textbf\let\PY@it=\textit}
\expandafter\def\csname PY@tok@kd\endcsname{\let\PY@bf=\textbf}
\expandafter\def\csname PY@tok@s1\endcsname{\let\PY@it=\textit}
\expandafter\def\csname PY@tok@sc\endcsname{\let\PY@it=\textit}
\expandafter\def\csname PY@tok@ni\endcsname{\let\PY@bf=\textbf}
\expandafter\def\csname PY@tok@err\endcsname{\def\PY@bc##1{\setlength{\fboxsep}{0pt}\fcolorbox[rgb]{1.00,0.00,0.00}{1,1,1}{\strut ##1}}}
\expandafter\def\csname PY@tok@ow\endcsname{\let\PY@bf=\textbf}
\expandafter\def\csname PY@tok@si\endcsname{\let\PY@bf=\textbf\let\PY@it=\textit}
\expandafter\def\csname PY@tok@sh\endcsname{\let\PY@it=\textit}
\expandafter\def\csname PY@tok@nt\endcsname{\let\PY@bf=\textbf}

\def\PYZus{\char`\_}
\def\PYZob{\char`\{}
\def\PYZcb{\char`\}}
\def\PYZca{\char`\^}

\def\PYZgt{\char`\>}

\def\PYZhy{\char`\-}

\makeatother
\makeatletter{}% solid tombstone QED symbols
\newcommand{\proofqedsymbol}{$\blacksquare$}
\newcommand{\closeqedsymbol}{$\blacklozenge$}

\renewcommand{\qedsymbol}{\proofqedsymbol}

\makeatletter
\newcommand{\setlabeltext}[1]{\def\@currentlabel{#1}}
\makeatother

\newcounter{theoremnum}[section]
\newcommand{\thetheorem}{\thesection.\arabic{theoremnum}}

\newlength{\mytemplength}

\newcommand{\mynewtheorem}[2]{
  \newenvironment{#1}[1][]{
    \stepcounter{theoremnum}
    \settowidth{\mytemplength}{##1}
    \renewcommand{\qedsymbol}{\closeqedsymbol}
    \ifthenelse{\lengthtest{\mytemplength=0pt}}{
      \begin{proof}[\normalfont\bfseries #2 \thetheorem]
    }{
      \begin{proof}[{\normalfont\textbf{#2 \thetheorem} ({##1})}\normalfont\bfseries]
    }
    \setlabeltext{\thetheorem}
  }{
    \end{proof}
    \renewcommand{\qedsymbol}{\proofqedsymbol}
  }
}

\mynewtheorem{theorem}{Theorem}
\mynewtheorem{lemma}{Lemma}
\mynewtheorem{proposition}{Proposition}
\mynewtheorem{corollary}{Corollary}
\mynewtheorem{claim}{Claim}
\mynewtheorem{definition}{Definition}
\mynewtheorem{example}{Example}
\mynewtheorem{remark}{Remark}
\newcommand{\FINAL}{}

\ifdefined\FINAL
\else
    
\fi

\makeatletter{}% proof macro
\newenvironment{pf}[1][\proofname]
{\proof[\normalfont \bfseries #1]}
{\endproof}

\newcommand{\lin}{\operatorname{lin}}

\newcommand{\eps}{\varepsilon}

\newcommand{\p}{\partial}
\renewcommand{\t}[1]{\widetilde{#1}}
\newcommand{\h}[1]{\hat{#1}}

\newcommand{\R}{\mathbf{R}}
\newcommand{\N}{\mathbf{N}}
\newcommand{\Z}{\mathbf{Z}}
\newcommand{\cN}{\mathcal{N}}
\newcommand{\M}{\mathcal{M}}

\newcommand{\T}{\mathcal{T}}

\newcommand{\dx}{\; dx}

\newcommand{\ds}{\; ds}

\newcommand{\pr}{\partial_r}
\newcommand{\px}{\partial_x}
\newcommand{\pt}{\partial_t}
\newcommand{\ps}{\partial_s}

\newcommand{\hx}{\hat{x}}
\newcommand{\he}{\hat{e}}

\newcommand{\ol}[1]{\overline{#1}}

\numberwithin{equation}{section}

\renewcommand{\qedsymbol}{$\blacksquare$}

\begin{document}
\title{Critical $O(d)$-equivariant biharmonic maps}
\author{Matthew K. Cooper}
\date{}
\makeatletter{}\begin{abstract}
We study $O(d)$-equivariant biharmonic maps in the critical dimension.
A major consequence of our study concerns the corresponding heat flow.
More precisely, we prove that blowup occurs in the biharmonic map heat flow from $B^4(0, 1)$ into $S^4$.
To our knowledge, this was the first example of blowup for the biharmonic map heat flow.
Such results have been hard to prove, due to the inapplicability of the maximum principle in the biharmonic case.
Furthermore, we classify the possible $O(4)$-equivariant biharmonic maps from $\R^4$ into $S^4$, and we show that there exists, in contrast to the harmonic map analogue, equivariant biharmonic maps from $B^4(0,1)$ into $S^4$ that wind around $S^4$ as many times as we wish.
We believe that the ideas developed herein could be useful in the study of other higher-order parabolic equations.

\medskip
\noindent\textbf{Mathematics Subject Classification}
\hspace{1ex}
35J40 $\cdot$
35J55 $\cdot$
35J60 $\cdot$
35J65 $\cdot$
58E20 $\cdot$
35K35 $\cdot$
35K55 $\cdot$
35B40 $\cdot$
34C11 $\cdot$
34C30
\end{abstract}
%%% Local Variables:
%%% mode: latex
%%% TeX-master: "main"
%%% End:

\maketitle

\section{Introduction}
\makeatletter{}In this work we study (extrinsic) biharmonic maps in the critical dimension.
Some of our results have applications to the corresponding heat flow.
Therefore, we first describe the biharmonic maps and their corresponding heat flow.

Let $\Omega \subset \R^d$ be open and bounded, and $\cN$ a smooth compact Riemannian manifold without boundary which is isometrically embedded in $\R^{\kappa}$ for some $\kappa \in \N$.
Consider the \emph{bi-energy}
\begin{equation}\label{BiEnergy}
E_2(u; \Omega) = \int_{\Omega} |\Delta u|^2 \dx
\text{ for }
u \in H^2(\Omega; \cN),
\end{equation}
where
\begin{equation*}
H^m(\Omega; \cN)
=
\{u \in H^m(\Omega; \R^{\kappa}) \; : \; u(x) \in \cN \text{ for a.e.\ } x \in \Omega \},
\end{equation*}
for $m \in \N_0$.

Critical points of this bi-energy are called biharmonic maps.
The Euler-Lagrange equation associated to \eqref{BiEnergy} is
\begin{equation}\label{BiEulerLagrangeEq}
\left\{
\begin{aligned}
&\Delta^2 u(x) \perp \T_{u(x)} \cN \mbox{ on } \Omega,
\\
&D^{\alpha} u = D^{\alpha} g \mbox{ on } \p \Omega \mbox{ for } |\alpha| \leq 1,
\end{aligned}
\right.
\end{equation}
where $g$ is boundary data, and the above is interpreted in the distributional sense for $u \in H^2(\Omega; \cN)$.

The $L^2$-gradient flow, or heat flow, of \eqref{BiEnergy} is
\begin{equation}\label{BMHFEq}
\left\{
\begin{aligned}
&\pt u(t,x) + \Delta^2 u(t,x)
\perp
\T_{u(t, x)} \cN \mbox{ on } \R^+ \times \Omega,
\\
&D^{\alpha} u = D^{\alpha} g \mbox{ on } \Gamma(\R^+ \times \Omega) \mbox{ for } |\alpha| \leq 1,
\end{aligned}
\right.
\end{equation}
where $g$ is initial-boundary data, and $\Gamma(U)$ denotes the parabolic boundary of $U \subset \R^{1+d}$.
We may replace $\Omega$ with a Riemannian manifold $\M$.
However, for concreteness, and since our study does not need this generality, our presentation will only consider flat domains.

The bi-energy \eqref{BiEnergy} is a higher-order analogue of the Dirichlet energy
\begin{equation*}
E_1(u; \Omega) = \int_{\Omega} |D u|^2 \dx,
\end{equation*}
for $u \in H^1(\Omega; \cN)$, where $Du$ is the Jacobian of $u$.
Therefore, one can view biharmonic maps as a higher-order analogue of harmonic maps.
Different higher-order energies have been proposed.
For example, the intrinsic bi-energy
\begin{equation}\label{IntrinsicBiEnergy}
H_2(u;\Omega) = \int_{\Omega} |(\Delta u(x))^T|^2 \dx,
\end{equation}
where $(\Delta u(x))^T$ is the orthogonal projection of $\Delta u(x)$ onto $\T_{u(x)} \cN$.
The intrinsic energy does not depend on the embedding of $\cN$ into Euclidean space, whereas the extrinsic energy, that is, the energy given by \eqref{BiEnergy}, does.
This makes the intrinsic energy more natural from a geometric perspective.
However, the intrinsic energy lacks coercivity, in contrast to the extrinsic energy, which makes it difficult to work with from an analytic perspective.
In this work we focus solely on the extrinsic energy.

Before we state our main results we need to introduce the notion of $O(d)$-equivariance that we use.
Throughout this paper $O(d)$ refers to the standard group of orthogonal transformations acting on $\R^d$.
We suppose that $u : \Omega \rightarrow S^d$, where $\Omega \subset \R^d$ is invariant under the action of $O(d)$, and $S^d$ is embeddeded in $\R^{d+1}$:
\begin{equation*}
S^d = \{x \in \R^{d+1} \; : \; |x| = 1\}.
\end{equation*}
For $x \in \R^d$ and $R \in O(d)$, we let $R x$ denote the standard group action of $O(d)$ on $\R^d$.
For $y = (\t{y}, y_{d+1}) \in \R^{d+1}$, we set $R \bullet y = (R \t{y}, y_{d+1})$.
We say that $u$ is $O(d)$-equivariant, or simply \emph{equivariant}, if $R \bullet u(x) = u(Rx)$, for all $R \in O(d)$ and $x \in \Omega$.

We will prove in Lemma \ref{OisS1Lemma} that, if $u \in C^\infty(\ol{B^d}; S^d)$ is equivariant then there exists a $k \in \N$ such that there exists a unique $\psi \in C^\infty([0, 1]; \R)$, called the longitudinal distance, where $\psi(0) = k \pi$ and
\begin{equation}\label{EquiAnsatzEq}
u(x)
=
\Upsilon(\psi)(x)
=
\begin{cases}
\left( \frac{x}{|x|} \sin\psi(|x|), \cos\psi(|x|) \right) &\mbox{ for } x \in \overline{B^d(0,1)} \setminus \{0\},
\\
\pm_k  \hat{e}_{d+1} &\mbox{ if } x = 0,
\end{cases}
\end{equation}
where
\begin{equation*}
\pm_k
=
\begin{cases}
+ &\text{ if } k \text{ is even},
\\
- &\text{ if } k \text{ is odd}.
\end{cases}
\end{equation*}

\begin{remark}
The condition $\psi(0) = k \pi$ ensures continuity of $u$ at the origin, if $\psi$ is itself continuous at zero.
\end{remark}

\begin{remark}\label{remark:PsiOriginCondition}
From \eqref{EquiAnsatzEq} it is clear that, if we set $\psi_0(r) = \psi(r) + 2 \pi l$, for $r \in [0, 1]$ and $l \in \Z$, then $\Upsilon(\psi_0) = \Upsilon(\psi)$.
Therefore, without loss of generality, we may assume that $\psi(0) \in \{0, \pi\}$.
\end{remark}

\begin{remark}\label{remark:PsiOriginConditionReflection}
Suppose that $\psi(0) = \pi$ and $u = \Upsilon(\psi)$.
If we set $\psi_0(r) = \pi - \psi(r)$ and $u_0 = \Upsilon(\psi_0)$ then $u_0(x) = \mathfrak{R}u(x)$, where $\mathfrak{R}$ is the reflection through the $\{x_{d+1} = 0\} \subset \R^{d+1}$ hyperplane.
\end{remark}

In this paper we globally assume that $\psi(0) = 0$.
This is without loss of any generality, because of remarks \ref{remark:PsiOriginCondition} and \ref{remark:PsiOriginConditionReflection}.
For clarity, we will occasionally remind the reader that we are assuming $\psi(0) = 0$.

Next, we describe our primary result.
For an equivariant map from $\ol{B^4(0, 1)}$ into $S^4$, if the normal derivative at the boundary vanishes then there is a limit on the number of times an equivariant biharmonic map from $\ol{B^4(0,1)}$ into $S^4$ satisfying the same boundary condition can wind around $S^4$.

\begin{theorem}\label{MainNonExistence}
There exists a $K > 0$ such that if
\begin{equation*}
u = \Upsilon(\psi) \in C^{\infty}(\ol{B^4(0,1)}; S^4),
\end{equation*}
with $\psi(0) = 0$, $|\psi(1) | \geq K$, and $\pr \psi(1) = 0$, then $u$ cannot be a biharmonic map.
\end{theorem}

This implies that the critical equivariant biharmonic map heat flow starting from such initial data must blowup in finite time or at infinity, since it cannot sub-converge to a biharmonic map.
To our knowledge, this was the first blowup result for the biharmonic map heat flow.

Recently in \cite{ChenLi13}, an example of finite-time blowup for the harmonic map heat flow due to topological reasons was given.
Their argument is based on a no-neck theorem, and builds upon earlier observations in \cite{QingTian97}.
In \cite{BiharmonicNeckAnalysis}, Liu and Yin prove a no-neck theorem for the blowup of a sequence of extrinsic and intrinsic biharmonic maps with bounded energy.
Motivated by the arguments in \cite{ChenLi13}, Liu and Yin, in \cite{BiharmonicFiniteTimeBlowup}, have used their no-neck theorem to show finite-time blowup for the biharmonic map heat flow in the critical dimension.
More precisely, they prove the following:

\begin{theorem}[{\cite[Theorem 1.1]{BiharmonicFiniteTimeBlowup}}]
Suppose that $\cN'$ is any closed manifold of dimension $n' > 4$ with nontrivial $\pi_4(\cN')$, and let $\cN = \cN' \;\#\; T^{n'}$ be the connected sum of $\cN'$ with the torus of the same dimension.
For any Riemannian metric $g$ on $\cN$, there exists (infinitely many) initial maps $u_0 : S^4 \rightarrow \cN$ such that the biharmonic map heat flow starting from $u_0$ develops a singularity in finite time.
\end{theorem}

It must be emphasized that the question of finite-time blowup for the biharmonic, and other higher-order polyharmonic, heat flows into spheres is still open.
It is in this latter case that we expect symmetric/equivariant solutions to play an important role.

With our second result we show that non-constant equivariant biharmonic maps from $\R^4$ into $S^4$ are unique, modulo dilation, reflection through the origin in the domain, and reflection through the plane $\{x_5 = 0\}$ in the codomain (see Remark \ref{remark:PsiOriginConditionReflection}).
More precisely, these transformations are, respectively, $u \mapsto (x \mapsto u(\lambda x))$, for $\lambda \in \R^+$, $u \mapsto (x \mapsto u(-x))$, and $u \mapsto (x \mapsto \mathfrak{R}u(x))$, where $\mathfrak{R}$ is the reflection through the plane $\{x_5 = 0\}$.

\begin{theorem}\label{MainUniqueBiharmMap}
If $u = \Upsilon(\psi) \in C^{\infty}(\R^4; S^4)$ is a non-constant equivariant biharmonic map with $\psi(0) = 0$ then, up to dilation, $\psi(r) = \pm 2 \arctan r$.
\end{theorem}

This is interesting, because in an upcoming paper the author shows that in the critical equivariant biharmonic map heat flow blowup coincides with a bubble separating at the origin.
Furthermore, this bubble is a non-constant equivariant biharmonic map from $\R^4$ into $S^4$, hence Theorem \ref{MainUniqueBiharmMap} gives a complete description of these bubbles.
By computing the energy of this bubble, we know that if our equivariant initial data has bi-energy less than or equal to $32\,\text{vol}(S^3)$ then the resulting flow exists globally in time and sub-converges to a smooth equivariant biharmonic map from $B^4(0,1)$ into $S^4$.

With our final result we show that, in contrast to the harmonic case, there are equivariant biharmonic maps from $B^4(0,1)$ into $S^4$ that wind around $S^4$ as many times as we wish.

\begin{theorem}\label{MainArbitraryWrapAround}
Let $a \in \R$.
Then there exists a biharmonic map
\begin{equation*}
u = \Upsilon(\psi) \in C^{\infty}(\ol{B^4(0,1)}; S^4),
\end{equation*}
such that $\psi(0) = 0$ and $\psi(1) = a$.
\end{theorem}

There is a growing volume of literature concerning biharmonic maps, and more generally higher-order polyharmonic maps, and their heat flows.
We will not attempt to survey this literature here, but the reader could start with
\cite{Angelsberg06},
\cite{CWY99},
\cite{FGOT12},
\cite{Gastel06},
\cite{GastelScheven09},
\cite{GastelZorn12},
\cite{GSZG09},
\cite{HHW2014},
\cite{Ku08},
\cite{LammSmallEnergy04},
\cite{Lamm05},
\cite{LammWang09},
\cite{BiharmonicFiniteTimeBlowup},
\cite{MR2013},
\cite{IntrinsicBlowupBehaviour},
\cite{WeakIntrinsic},
\cite{PolyharmonicUniquenessReverseBubbling},
\cite{Strzelecki03},
\cite{Wang04Arbitrary},
\cite{Wang04Spheres},
\cite{Wang04Stationary},
\cite{Chang07},
\cite{BiharmonicRough},
\cite{WOY2014},
\cite{ZornThesis},
and the references therein.

Symmetric and equivariant biharmonic maps have already been studied in \cite{GastelZorn12,MR2013,WOY2014,ZornThesis}.
In \cite{WOY2014}, Wang, Ou, and Yang study rotationally symmetric intrinsic biharmonic maps from $S^2$ into $S^2$.
Similar to this work, they compute the corresponding symmetry reduction, and classify their class of symmetric intrinsic biharmonic maps.
In \cite{MR2013}, Montaldo and Ratto examine a more general class of equivariant intrinsic biharmonic maps.
They consider maps that are equivariant with respect to Riemannian submersions.
They setup machinery to compute the corresponding symmetry reductions, and use this to explicitly compute the symmetry reduction in some concrete cases.
As applications they prove the stability of specific proper, that is, non-harmonic, intrinsic biharmonic maps from $T^2$ into $S^2$ among a certain class of equivariant maps.
Moreover, they construct a counter-example to a generalization to intrinsic biharmonic maps of Sampson's maximum principle for harmonic maps, see \cite{Sampson78}.
In \cite{ZornThesis}, Zorn studies $G$-equivariant biharmonic maps, where $G$ is a compact Lie group.
Among other results, he proves that $G$-minimizers of the bi-energy are stationary biharmonic, and improves estimates on the Hausdorff dimension of the singular sets of appropriate $G$-equivariant biharmonic maps.

Biharmonic maps are related to interesting problems in four dimensional conformal geometry, for example see \cite{ChangGurskyYang99} and \cite{XuYang02}.
A further motivation for the study of biharmonic maps and their heat flows is that they are good model equations for other interesting higher-order elliptic and parabolic PDE.
Therefore, one would hope that the insights generated for biharmonic maps and their heat flows to be useful elsewhere.
As an example of the parallels between the biharmonic map, and more generally polyharmonic map, heat flow and other higher-order parabolic equations, one may compare the work on the Willmore flow by Kuwert and Sch\"atzle, see \cite{KuwertSchatzle01}, and the work on the polyharmonic map heat flow and biharmonic map heat flow respectively by Gastel, see \cite{Gastel06}, and Wang, see \cite{Chang07}.

In \cite{ChangDingYe1992} Chang, Ding, and Ye showed that finite-time blowup was possible for the critical harmonic map heat flow.
Chang, Ding, and Ye worked with the equivariant ansatz \eqref{EquiAnsatzEq}.
More precisely, Chang, Ding, and Ye's result says that if $\psi(0, 0) = 0$ and $|\psi(0,1)| > \pi$ then the corresponding solution to the critical harmonic map heat flow blows up in finite time.
Prior to this, Chang and Ding, in \cite{ChangDing91}, show that if $\psi(0, 0) = 0$ and $\sup_{r \in [0,1]} |\psi(0,r)| \leq \pi$ then the corresponding local in time classical solution to the critical harmonic map heat flow is in fact global in time.
Only when $|\psi(0,1)| < \pi$ do we have sub-convergence to a harmonic map, see \cite{KarcherWood84} for further details.
These proofs rely heavily on the comparison/maximum principle.

Our main motivation for this study is to try to extend the work in \cite{ChangDing91,ChangDingYe1992} to the biharmonic case.
The inapplicability of the maximum principle in the biharmonic case makes this task difficult.
Looking to the future the following ideas may be useful.
In \cite{GalPoh02} Galaktionov and Pohozaev use the technique of majorizing operators in order to obtain a comparison principle for the bi-heat equation.
These ideas may allow us to construct barriers in order to prove global existence in a similar way as in \cite{ChangDing91}.
However, this will not allow a proof of finite-time blowup using barriers.
In \cite{PierreSchweyer13} Rapha{\"e}l and Schweyer look at finite-time blowup of the $1$-corotational critical harmonic map heat flow while avoiding the maximum principle.
Instead they rely on energy methods and modulation theory.

Recently in \cite{GastelZorn12} Gastel and Zorn studied a fourth-order ODE very similar to our symmetry reduction for the equivariant biharmonic map equation \eqref{BMHFSymRedEq} (with $\pt \psi = 0$).
Their ODE arises when trying to construct biharmonic maps of cohomonogeneity one between spheres using joins of two harmonic eigenmaps.
In contrast to this work, they use purely variational methods.
Their stated reason for this choice being that purely ODE methods would cause difficulties due to their ODE being fourth-order and having ill-posed boundary conditions.
Although the questions studied in \cite{GastelZorn12} are quite different than the ones studied here, we have found success in using ODE methods to study our similar fourth-order ODE.
To us, this demonstrates that ODE methods can be useful in exploring such questions.
It may be of interest in future work to see if a synthesis of the ideas in \cite{GastelZorn12} and here can yield deeper insights.

We will now outline the structure of the rest of this paper.
In Section \ref{EquivariantSection}, we first prove that the smooth flow of \eqref{BMHFEq} preserves $O(d)$-equivariance.
In \cite[Lemma 2.2]{ChangDing91}, the authors study the harmonic map heat flow using the equivariant ansatz \eqref{EquiAnsatzEq}, and they prove the analogue of this result in their situation using the maximum principle, see also \cite[Lemma 4.2]{Grotowski91} where a similar argument using the maximum principle is used for the axially symmetric harmonic map heat flow.
Unfortunately, the maximum principle is not available in the biharmonic map heat flow.
Therefore, in our setting a different method which avoids the maximum principle must be used.
Next, we show the equivalence of maps which are $O(d)$-equivariant and those which satisfy \eqref{EquiAnsatzEq}.
After this, we present Mathematica code for computing our symmetry reduction of the biharmonic map heat flow.
This is then used to explicitly compute our symmetry reduction in the critical dimension.

In Section \ref{EquiBiharmonicSection}, we prepare for, and outline our approach to, our deeper study of equivariant biharmonic maps in the critical dimension.
This study is carried out in sections \ref{BlowupOrHetroSection} and \ref{UnstableManifoldSect}, and it is here that the theorems mentioned above are proven.
The paper ends with an appendix which contains a couple of technical proofs which, in the author's opinion, do not yield much conceptual insight.

{\bf Notation.}
Throughout this paper $C$ denotes a positive universal constant.
Two different occurrences of $C$ are liable to be different.
If our constant depends on some parameter, say $\eps$, then we denote this by writing $C(\eps)$.

{\bf Acknowledgments.}
This work was undertaken while the author was a PhD student at The University of Queensland, Australia under the supervision of Prof. Joseph Grotowski and Prof. Peter Adams.
The author owes much of his success to the advice, guidance, and support of his supervisors.

This work was supported by an Australian Postgraduate Award and the Australian Research Council (discovery grant DP120101886).
The author would like to thank Prof. Yihong Du, Prof. Dr. Andreas Gastel, Prof. Joseph Grotowski, and the anonymous referees for comments and suggestions that have improved this paper.
%%% Local Variables:
%%% mode: latex
%%% TeX-master: "main"
%%% End:

\section{The equivariant ansatz}\label{EquivariantSection}
\makeatletter{}In this section we show that $O(d)$-equivariance is preserved by the smooth biharmonic map heat flow, and $O(d)$-equivariance is equivalent to the ansatz given by \eqref{EquiAnsatzEq}.
After this, Mathematica code for computing the symmetry reduction is presented, and then used to explicitly compute the symmetry reduction for the biharmonic map heat flow in the critical dimension.
%%% Local Variables:
%%% mode: latex
%%% TeX-master: "../main"
%%% End:

\makeatletter{}
First we show that $O(d)$-equivariance is preserved by the smooth biharmonic map heat flow.

\begin{lemma}\label{GInvPreservedLem}
Suppose that $d \in \N_{\geq 2}$, $T > 0$, $Q = [0,T] \times \ol{B^d(0,1)}$, and $u \in C^{\infty}(Q; S^d)$ is a solution to \eqref{BMHFEq} with $g$ as $O(d)$-equivariant initial-boundary data.
Then, for each $t \in [0,T]$, $u(t,\cdot)$ is $O(d)$-equivariant.
\end{lemma}

\begin{pf}
We let $R \in O(d)$ be arbitrary, and set
\begin{equation*}
v_R(t,x) = R \bullet u(t, R^{-1} x) \text{ for } (t,x) \in Q.
\end{equation*}
From \eqref{BMHFEq}, we have
\begin{equation*}
\pt u = -\Delta^2 u + (\Delta^2 u \cdot u) u,
\end{equation*}
hence
\begin{equation*}
\pt v_R(t,x)
=
-\Delta^2 v_R(t,r) + (\Delta^2 v_R(t,r) \cdot v_R(t,r)) v_R(t,r),
\end{equation*}
where we have used properties of orthogonal matrices.
Since $g$ is $O(d)$-equivariant, we have $D^{\alpha} v_R = D^{\alpha} g$ on $\Gamma Q$ for $|\alpha| \leq 1$.
Therefore, $v_R$ solves \eqref{BMHFEq} with the same initial-boundary data.
Since we are working in the smooth category, uniqueness of solutions is standard, hence $v_R \equiv u$ and $u(t,x) = R \bullet u(t, R^{-1} x)$ for all $(t,x) \in Q$ and $R \in O(d)$.
\end{pf}

\begin{remark}
With the obvious modifications the above proof works for polyharmonic maps of any order.
\end{remark}

Next, we have two lemmas that demonstrate the equivalence between $O(d)$-equivariance and the ansatz given by \eqref{EquiAnsatzEq}.

\begin{lemma}\label{S1EquiImpliesGInvLem}
Suppose that $k \in \Z$, $d \in \N_{\geq 2}$, $u \in C(\ol{B^d(0,1)}; S^d)$ satisfies \eqref{EquiAnsatzEq}.
Then $u$ is $O(d)$-equivariant.
\end{lemma}

\begin{pf}
We let $R \in O(d)$ be arbitrary.
Using \eqref{EquiAnsatzEq}, we compute that
\begin{equation*}
R \bullet u(R^{-1} x) = u(x) \text{ for } x \neq 0.
\end{equation*}
For $x = 0$, we have
\begin{equation*}
R \bullet u(R^{-1} 0) = R \bullet u(0) = R \bullet \left(\pm_k \he_{d+1}\right) = \pm_k \he_{d+1} = u(0).
\qedhere
\end{equation*}
\end{pf}

Finally, we show that if a map is $O(d)$-equivariant at each time then it satisfies \eqref{EquiAnsatzEq} at each time.

\begin{lemma}\label{OisS1Lemma}
Let $d \in \N_{\geq 2}$, $k \in \Z$, and $T > 0$.
Suppose that $u \in C^{\infty}([0,T] \times \ol{B^d(0,1)}; S^d)$, $u(t,\cdot)$ is $O(d)$-equivariant and $u(t,0) = \pm_k \he_{d+1}$ for each $t \in [0,T]$.
Then there exists a unique $\psi \in C^{\infty}([0,T] \times [0,1])$ such that $\psi(t,0) = k \pi$ and $u(t, \cdot) = \Upsilon(\psi(t, \cdot))$ for each $t \in [0,T]$.
\end{lemma}

\begin{pf}
Let $R_0 \in O(d)$ be the reflection through the $\he_1$ axis, that is,
\begin{equation}\label{ReflectThroughBasisVecEq}
R_0 \he_i
=
\begin{cases}
\he_1 &\mbox{ if } i = 1,
\\
-\he_i &\mbox{ otherwise}.
\end{cases}
\end{equation}
For $r \in [0,1]$, we have
\begin{equation}\label{OisS1Eq1}
u(t,r \he_1) = R_0 \bullet u(t,R_0^{-1} r \he_1) = R_0 \bullet u(t,r\he_1).
\end{equation}
Therefore, $u^i(t,r \he_1) = 0$ for $i \in \{2,3,\dots,d\}$, and
\begin{equation*}
(u^1(t,r \he_1),u^{d+1}(t,r \he_1)) \in S^1.
\end{equation*}
Since $u \in C^{\infty}([0,T] \times \overline{B^d(0,1)}; S^d)$ and $u(t,0)=\pm_k \he_{d+1}$, there exists a unique $\psi \in C^{\infty}([0,T] \times [0,1])$ such that $\psi(t,0) = k \pi$ and
\begin{equation}\label{OisS1Eq2}
u^1(t,r \he_1) = \sin \psi(t,r), \quad u^{d+1}(t,r \he_1) = \cos \psi(t,r).
\end{equation}
Next, we work in spherical coordinates.
We fix an $\hx \in S^{d-1}$, and let $R_{\hx} \in O(d)$ be a map such that $R_{\hx} \he_1 = \hx$.
Then for $r > 0$ we calculate, keeping in mind \eqref{OisS1Eq1} and \eqref{OisS1Eq2},
\begin{align*}
u(t,r\hx)
&=
R_{\hx} \bullet u(t,r R_{\hx}^{-1} \hx)
\\
&=
R_{\hx} \bullet  u(t,r \he_1)
\\
&=
R_{\hx} \bullet ( \sin \psi(t,r) \; \he_1 + \cos \psi(t,r) \; \he_{d+1} )
\\
&=
\Upsilon(\psi(t,\cdot))(r\hx).
\end{align*}
Since $\hx$ was arbitrary, we are done.
\end{pf}
%%% Local Variables:
%%% mode: latex
%%% TeX-master: "../main"
%%% End:

\makeatletter{}Next, we compute the symmetry reduction.
For a function $f : \R^d \rightarrow \R$ such that $f(x) = f(|x|) = f(r)$, we have
\begin{equation*}
\Delta f(x) = \pr^2 f(r) + \frac{d-1}{r} \pr f(r) =: (L_1 f)(r) \mbox{ for } x \neq 0.
\end{equation*}
We also compute:
\begin{equation*}
\Delta \left( \frac{x}{|x|} f(x) \right)
=
\frac{x}{|x|} \left( \pr^2 f(r) + \frac{d-1}{r} \pr f(r) - \frac{d-1}{r^2} f(r) \right)
=:
\frac{x}{|x|}  (L_0 f)(r)
\text{ for } x \neq 0.
\end{equation*}
For appropriate $g_0, g_1 : [0,1] \rightarrow \R$, we write
\begin{equation}\label{BraceNotationEq}
\{g_0, g_1\}_{\h{x}}(t,x)
=
(\h{x} g_0(t,|x|), g_1(t,|x|)).
\end{equation}
We observe that, for $x \neq 0$,
\begin{align*}
\Delta \{g_0, g_1\}_{\h{x}}
&=
\{L_0 g_0, L_1 g_1\}_{\h{x}},
\\
\pt \{g_0, g_1\}_{\h{x}}
&=
\{\pt g_0, \pt g_1\}_{\h{x}}, \mbox{ and }
\\
\{g_0, g_1\}_{\h{x}} \cdot \{h_0, h_1\}_{\h{x}}
&=
g_0 h_0 + g_1 h_1.
\end{align*}
For $\cN$ a unit sphere, \eqref{BMHFEq} can be written as
\begin{equation*}
|\pt u + \Delta^2 u - \left(\Delta^2 u \cdot u\right) u|^2 = 0.
\end{equation*}
When $u$ is $O(d)$-equivariant this reduces to a PDE for $\psi$.

The following Mathematica code computes this symmetry reduction in the critical dimension.
\begin{Verbatim}[commandchars=\\\{\},numbers=left,firstnumber=1,stepnumber=1,frame=lines,codes={\catcode`\$=3\catcode`\^=7}]
\PY{c}{(* Functions to calculate the Laplacian *)}
\PY{n}{L1}\PY{p}{[}\PY{n+nv}{expr\PYZus{}}\PY{p}{]}\PY{+w}{ }\PY{o}{:=}\PY{+w}{ }\PY{n}{D}\PY{p}{[}\PY{n}{expr}\PY{p}{,}\PY{+w}{ }\PY{p}{\PYZob{}}\PY{n}{r}\PY{p}{,}\PY{+w}{ }\PY{l+m+mi}{2}\PY{p}{\PYZcb{}}\PY{p}{]}\PY{+w}{  }\PY{o}{+}\PY{+w}{ }\PY{p}{(}\PY{n}{d}\PY{l+m+mi}{\PYZhy{}1}\PY{p}{)}\PY{o}{/}\PY{n}{r}\PY{+w}{ }\PY{n}{D}\PY{p}{[}\PY{n}{expr}\PY{p}{,}\PY{+w}{ }\PY{n}{r}\PY{p}{]}
\PY{n}{L0}\PY{p}{[}\PY{n+nv}{expr\PYZus{}}\PY{p}{]}\PY{+w}{ }\PY{o}{:=}\PY{+w}{ }\PY{n}{L1}\PY{p}{[}\PY{n}{expr}\PY{p}{]}\PY{+w}{ }\PY{o}{\PYZhy{}}\PY{+w}{ }\PY{p}{(}\PY{n}{d}\PY{l+m+mi}{\PYZhy{}1}\PY{p}{)}\PY{o}{/}\PY{n}{r}\PY{o}{\PYZca{}}\PY{l+m+mi}{2}\PY{+w}{ }\PY{n}{expr}
\PY{n}{Lapl}\PY{p}{[}\PY{n+nv}{expr\PYZus{}}\PY{p}{]}\PY{+w}{ }\PY{o}{:=}\PY{+w}{ }\PY{p}{\PYZob{}}\PY{n}{L0}\PY{p}{[}\PY{n}{expr}\PY{p}{[}\PY{p}{[}\PY{l+m+mi}{1}\PY{p}{]}\PY{p}{]}\PY{p}{]}\PY{p}{,}\PY{+w}{ }\PY{n}{L1}\PY{p}{[}\PY{n}{expr}\PY{p}{[}\PY{p}{[}\PY{l+m+mi}{2}\PY{p}{]}\PY{p}{]}\PY{p}{]}\PY{p}{\PYZcb{}}

\PY{c}{(* Calculate the symmetry reduction in the critical dimension *)}
\PY{n}{CriticalSymmReduction}\PY{+w}{ }\PY{o}{=}
\PY{+w}{  }\PY{n}{FullSimplify}\PY{p}{[}
\PY{+w}{    }\PY{n}{ReplaceAll}\PY{p}{[}
\PY{+w}{      }\PY{n}{With}\PY{p}{[}\PY{p}{\PYZob{}}\PY{n}{u}\PY{+w}{ }\PY{o}{=}\PY{+w}{ }\PY{p}{\PYZob{}}\PY{n}{Sin}\PY{p}{[}$\psi$\PY{p}{[}\PY{n}{t}\PY{p}{,}\PY{+w}{ }\PY{n}{r}\PY{p}{]}\PY{p}{]}\PY{p}{,}\PY{+w}{ }\PY{n}{Cos}\PY{p}{[}$\psi$\PY{p}{[}\PY{n}{t}\PY{p}{,}\PY{n}{r}\PY{p}{]}\PY{p}{]}\PY{p}{\PYZcb{}}\PY{p}{\PYZcb{}}\PY{p}{,}
\PY{+w}{        }\PY{n}{With}\PY{p}{[}\PY{p}{\PYZob{}}\PY{n}{v}\PY{+w}{ }\PY{o}{=}\PY{+w}{ }\PY{n}{Nest}\PY{p}{[}\PY{n}{Lapl}\PY{p}{,}\PY{+w}{ }\PY{n}{u}\PY{p}{,}\PY{+w}{ }\PY{l+m+mi}{2}\PY{p}{]}\PY{p}{\PYZcb{}}\PY{p}{,}
\PY{+w}{          }\PY{n}{With}\PY{p}{[}\PY{p}{\PYZob{}}\PY{n}{expr}\PY{+w}{ }\PY{o}{=}\PY{+w}{ }\PY{n}{D}\PY{p}{[}\PY{n}{u}\PY{p}{,}\PY{+w}{ }\PY{n}{t}\PY{p}{]}\PY{+w}{ }\PY{o}{+}\PY{+w}{ }\PY{n}{v}\PY{+w}{ }\PY{o}{\PYZhy{}}\PY{+w}{ }\PY{p}{(}\PY{n}{v}\PY{+w}{ }\PY{}{.}\PY{+w}{ }\PY{n}{u}\PY{p}{)}\PY{+w}{ }\PY{n}{u}\PY{p}{\PYZcb{}}\PY{p}{,}
\PY{+w}{            }\PY{n}{expr}\PY{+w}{ }\PY{}{.}\PY{+w}{ }\PY{n}{expr}\PY{+w}{ }\PY{o}{=}\PY{o}{=}\PY{+w}{ }\PY{l+m+mi}{0}\PY{p}{]}\PY{+w}{ }\PY{p}{]}\PY{p}{]}\PY{p}{,}
\PY{+w}{      }\PY{n}{d}\PY{+w}{ }\PY{o}{\PYZhy{}\PYZgt{}}\PY{+w}{ }\PY{l+m+mi}{4}\PY{p}{]}\PY{p}{]}
\end{Verbatim}
Exporting the result of this computation to \LaTeX{}, and rearranging, we obtain
\begin{equation}\label{BMHFSymRedEq}
\begin{aligned}
\pt \psi
&=
-
\pr^4 \psi
-
\frac{6}{r} \pr^3 \psi
+
6 (\pr \psi)^2 \pr^2 \psi
+
\frac{3}{r^2} \cos(2 \psi) \pr^2 \psi
+
\frac{6}{r} (\pr \psi)^3
\\
&\quad
-
\frac{3}{r^2} \sin(2 \psi) (\pr \psi)^2
+
\frac{3}{r^3}
\left(
\cos(2 \psi)
+
2
\right)
\p_r \psi
-
\frac{9}{2 r^4} \sin(2 \psi).
\end{aligned}
\end{equation}
Of course, we could have calculated \eqref{BMHFSymRedEq} by hand, but we prefer to delegate repetitive and elementary calculations to the computer.

Recall that we are assuming $\psi(0) = 0$, and hence $u(0) = \hat{e}_{d+1}$, see remarks \ref{remark:PsiOriginCondition} and \ref{remark:PsiOriginConditionReflection}.
Due to the boundary conditions in \eqref{BMHFEq}, we have
\begin{equation*}
\pr^i \psi(t,1) = a_i,
\end{equation*}
for $t \in [0,T]$ and $i \in \{0, 1\}$, where $a_i \in \R$.

Due to symmetry, $\psi$ must satisfy conditions at the origin.
Let $R_0$ be the same as in \eqref{ReflectThroughBasisVecEq}.
Arguing similarly as in Lemma \ref{OisS1Lemma}, we see that there exists a $\xi \in C^{\infty}([0,T] \times [-1,1];\R)$ such that $\xi(t, 0) = 0$ and
\begin{equation*}
u(t, x_1 \h{e}_1) = \h{e}_1 \sin \xi(t,x_1) + \h{e}_{d+1} \cos \xi(t,x_1),
\end{equation*}
for $x_1 \in [-1,1]$.
We also have
\begin{equation*}
u(t, x_1 \h{e}_1)
=
(-R_0) \bullet u(t, -x_1 \h{e}_1).
\end{equation*}
This implies that
\begin{equation*}
(\sin(\xi(t, x_1)), \cos(\xi(t, x_1)))
=
(\sin(-\xi(t, -x_1)), \cos(-\xi(t, -x_1))),
\end{equation*}
hence $\xi(t, \cdot)$ is odd.

Observe that if $u \in C^k([0, T] \times \ol{B^d(0,1)};S^d)$ then $\xi \in C^k([0, T] \times [-1,1];\R)$.
In this case, we have $\p_{x_1}^{2i} \xi(t, 0) = 0$ whenever $2i \leq k$.
Since $\psi = \xi|_{x \in [0,1]}$, $\psi \in C^k([0, T] \times [0,1];\R)$ and $\pr^{2i} \psi(t, 0) = 0$ whenever $2i \leq k$.
%%% Local Variables:
%%% mode: latex
%%% TeX-master: "../main"
%%% End:

\section{Critical equivariant biharmonic maps}\label{EquiBiharmonicSection}
\makeatletter{}Now we start to delve deeper into $O(d)$-equivariant biharmonic maps in the critical dimension.

After setting $\pt \psi = 0$ in \eqref{BMHFSymRedEq} and making the change of variables $\psi(r) = \phi(s(r))$, where $s(r) = \log r$, \eqref{BMHFSymRedEq} becomes the the fourth-order autonomous ODE
\begin{equation}\label{FourthOrderODE}
\ps^4 \phi
=
-\frac{9}{2} \sin(2 \phi) + (7 + 3 \cos(2\phi)) \ps^2 \phi
+
3 (\ps \phi)^2 (2 \ps^2 \phi - \sin(2\phi)).
\end{equation}
Recall that we assume $\psi(0) =0$.
This boundary condition becomes
\begin{equation}\label{TransformedBdryCond}
\lim_{s \rightarrow -\infty} \phi(s) = 0.
\end{equation}
We rewrite this as a first-order system by setting $\Phi_i = \ps^{i-1} \phi$ for $i \in \{1,2,3,4\}$:
\begin{equation}\label{FirstOrderODE}
\ps \Phi
=
\begin{pmatrix}
\Phi_2 \\
\Phi_3 \\
\Phi_4 \\
F_1(\Phi_1,\Phi_3)
+
\Phi_2^2 F_2(\Phi_1, \Phi_3)
\end{pmatrix},
\end{equation}
where
\begin{align*}
F_1(\Phi_1, \Phi_3) &= -\frac{9}{2} \sin(2 \Phi_1) + (7 + 3 \cos(2\Phi_1)) \Phi_3, \mbox{ and}
\\
F_2(\Phi_1, \Phi_3) &= 3 (2 \Phi_3 - \sin(2\Phi_1)).
\end{align*}
The boundary condition \eqref{TransformedBdryCond} becomes
\begin{equation}\label{TransformedBdryCond2}
\lim_{s \rightarrow -\infty} \Phi_1(s) = 0.
\end{equation}
Observe that \eqref{FirstOrderODE} and \eqref{TransformedBdryCond2} are invariant under the transformation $\Phi \mapsto - \Phi$.

An $s_0 \in \R$ and initial data $\Phi^0(s_0) \in \R^4$ generate a unique solution to \eqref{FirstOrderODE}, denoted by $\Phi^0 : [s_0, s_{\max}) \rightarrow \R^4$, where either $s_{\max} = \infty$ or $\lim_{s \nearrow s_{\max}} |\Phi^0(s)| = \infty$.
The next lemma shows the equivalence between a solution of \eqref{FirstOrderODE} satisfying \eqref{TransformedBdryCond2} and it being an orbit in the unstable manifold of the origin of \eqref{FirstOrderODE}, denoted from now on by $W^u(0)$.

\begin{lemma}\label{BdryCondUnstableLem}
Suppose that $\Phi^0 : (-\infty, s_{\max}) \rightarrow \R^4$ solves \eqref{FirstOrderODE}.
Then the following are equivalent:
\begin{enumerate}
\item
$\displaystyle\lim_{s \rightarrow -\infty} \Phi^0_1(s) = 0$;
\item
$\displaystyle\lim_{s \rightarrow -\infty} \Phi^0(s) = 0$.
\qedhere
\end{enumerate}
\end{lemma}

For the proof see the Appendix.

Observe that $y(s) = 2 \arctan(e^s)$ is a heteroclinic orbit of \eqref{FourthOrderODE}, and $\Upsilon(2 \arctan(\cdot))$, for $d \geq 2$, is the inverse of the stereographic projection of $S^{d-1} \setminus \{-\h{e}_d\}$ onto $\R^{d-1}$.
We set $Y_i = \ps^{i-1} y$ for $i \in \{1,2,3,4\}$.

The first result we will focus on concerning \eqref{FirstOrderODE} states that this heteroclinic orbit gives rise to the only non-constant equivariant biharmonic map from $\R^4$ into $S^4$ modulo dilation, reflection through the origin in the domain, and reflection through the plane $\{x_5 = 0\}$ in the codomain (see Remark \ref{remark:PsiOriginConditionReflection}).

\begin{theorem}\label{MainUniqueThm}
Suppose that $\Phi^0 : (-\infty, s_{\max}) \rightarrow \R^4$ is a non-trivial orbit in $W^u(0)$.
Then the following dichotomy holds:
\begin{enumerate}
\item
up to $s$-translation $\Phi^0(s) = Y(s)$ or $\Phi^0(s) = -Y(s)$, hence $s_{\max} = \infty$; or
\item
$\Phi^0$ blows up in finite time, that is, $s_{\max} < \infty$ and
\begin{equation*}
\lim_{s \nearrow s_{\max}} |\Phi^0(s)| = \infty.
\qedhere
\end{equation*}
\end{enumerate}
\end{theorem}

Theorem \ref{MainUniqueBiharmMap} is a corollary of this.
Next, we outline our strategy for the study of \eqref{FirstOrderODE}.

\subsection*{Strategy}
Our arguments are in part motivated by the harmonic map case.
The analogue of \eqref{FourthOrderODE} in the critical harmonic map case is
\begin{align}\label{HarmonicSecondOrderODE}
\ps^2 \phi_h
=
\frac{1}{2} \sin(2 \phi_h)
&=
-\partial_{\phi_h} \left(\frac{1}{2} \cos^2 \phi_h \right),
\\
\nonumber
\lim_{s \rightarrow -\infty} \phi_h(s) &= 0.
\end{align}
This ODE is a pendulum equation.
One can think of it as describing the dynamics of a ball rolling without friction in coordinate space on the potential energy surface $V(q) = \frac{1}{2} \cos^2 q$.
After some consideration it is clear that if $\lim_{s \rightarrow -\infty} \phi_h(s) = 0$ and $\phi_h$ is non-constant then up to $s$-translation there are only two possibilities for $\phi_h$.
These possibilities being the heteroclinic orbits between $(\phi_h, \ps \phi_h) = (0, 0)$ and $(\phi_h, \ps \phi_h) = (\pi, 0)$, and between $(\phi_h, \ps \phi_h) = (0, 0)$ and $(\phi_h, \ps \phi_h) = (-\pi, 0)$.
These orbits happen to be, up to $s$-translation, $\pm y(s)$.

We see that if $\phi_h$ is an $s$-translation of $\pm y(s)$ then $|\phi_h| < \pi$.
Therefore, if we have equivariant initial data $u_0 = \Upsilon(\psi_0)$ for the harmonic map heat flow from $B^2(0,1)$ into $S^2$ such that $\psi_0(0) = 0$ and $|\psi_0(1)| \geq \pi$ then the flow must blowup either in finite time or at infinity, because it cannot sub-converge to a harmonic map.
Theorem \ref{MainNonExistence} is the analogue of this observation for the biharmonic map case.

\begin{figure}[t!]
\centering
\includegraphics[width=0.75\textwidth]{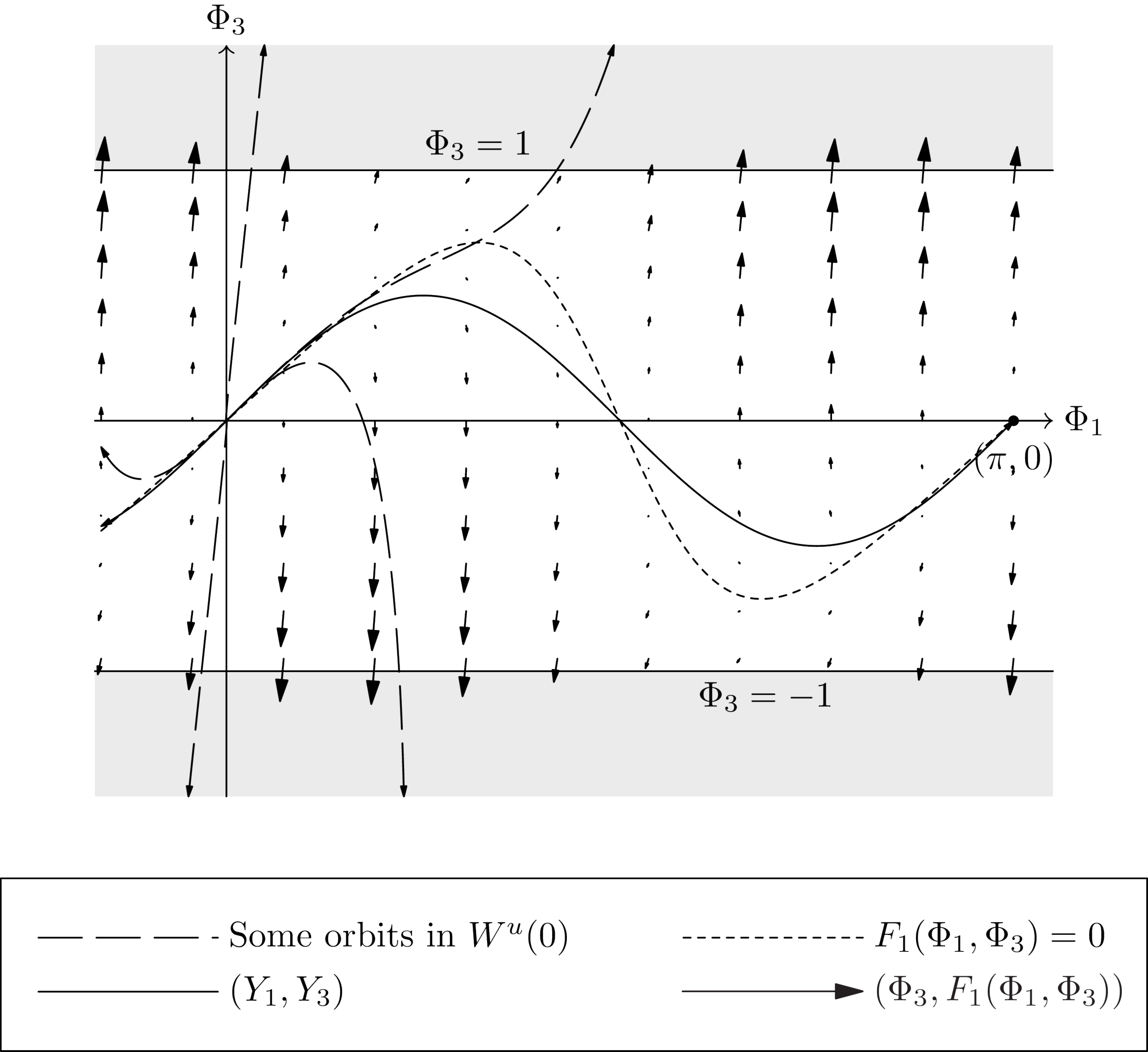}
\caption{%
This represents the portion of the coordinate space that is of interest for mid life orbits in $W^u(0)$.
We do not display the whole left half of the region, because this can be deduced, via symmetry, from the right half.
Note that there is an additional force acting in the vertical direction with a magnitude of $(\ps \Phi_1)^2 F_2(\Phi_1, \Phi_3)$.
This force acts to repulse solutions away from the graph of $\Phi_3 = \frac{1}{2} \sin(2 \Phi_1)$ in the coordinate space.
Observe that $Y_3 = \frac{1}{2} \sin(2 Y_1)$.
}
\label{figure:CoordinateSpace}
\end{figure}

Our situation is more complicated than the one encountered when studying \eqref{HarmonicSecondOrderODE}, because instead of a one-dimensional coordinate space we now have a two-dimensional coordinate space.
This adds complexity to the possible dynamics.
Moreover, unlike \eqref{HarmonicSecondOrderODE}, the dynamics of \eqref{FourthOrderODE} seem not to be related to a simple dynamical system from which we can gain intuition.
In spite of this, the author found it fruitful to consider \eqref{FourthOrderODE} as the following coupled system of second-order ODE:
\begin{equation*}
\left\{
\begin{aligned}
\ps^2 \Phi_1 &= \Phi_3
\\
\ps^2 \Phi_3 &=
F_1(\Phi_1, \Phi_3)
+
(\ps \Phi_1)^2 F_2(\Phi_1, \Phi_3),
\end{aligned}
\right.
\end{equation*}
and to think about the `forces' acting on the system in the coordinate space $(\Phi_1, \Phi_3)$, see Figure \ref{figure:CoordinateSpace}.

Our arguments are also inspired by the ideas in \cite{Gastel02}.
There the shooting method along with a pendulum equation interpretation was used to show the existence of singularities of the first kind in the harmonic map and Yang-Mills heat flows.

Much of our analysis revolves around finding positive invariant sets on which we approximate \eqref{FirstOrderODE} by systems of simpler ordinary differential inequalities that still give us enough control over the orbits in $W^u(0)$.
We found it convenient to divide the life of each orbit in $W^u(0)$ into three stages:

{\bf Early life.}
This is when the orbit is still close to the origin and its dynamics are well approximated by the linearization of \eqref{FirstOrderODE} at the origin.

{\bf Mid life.}
This is the most delicate stage to analyze, because when $|\Phi_3|$ is not so large we must deal with difficulties caused by the variation of the `forces' acting on $\Phi$ with respect to $\Phi_1$, see Figure \ref{figure:CoordinateSpace}.

{\bf Late life.}
This is when $|\Phi_3| \geq 1$.
In this case the variation of the `forces' with respect to $\Phi_1$ can be ignored.
This simplifies the situation, and one can prove finite-time blowup of an orbit in $W^u(0)$ once it reaches late life.
%%% Local Variables:
%%% mode: latex
%%% TeX-master: "main"
%%% End:

\section{Finite-time blowup or heteroclinic orbit}\label{BlowupOrHetroSection}
\makeatletter{}The aim of this section is to prove Theorem \ref{MainUniqueThm}, and to collect some facts along the way which will be used in our later arguments.
We let
\begin{align*}
W_+ &= \{ x \in \R^4 \; : \; (x_1,x_3), (x_2, x_4) \in \Lambda_+\} \mbox{ and}
\\
W_- &= \{ x \in \R^4 \; : \; (x_1,x_3), (x_2, x_4) \in \Lambda_-\},
\end{align*}
where
\begin{align*}
\Lambda_+ &= \{ x \in \R^2 \; : \; x_1 \geq 0, x_2 \geq 2 x_1\} \mbox{ and}
\\
\Lambda_- &= \{ x \in \R^2 \; : \; x_1 \leq 0, x_2 \leq 2 x_1\}.
\end{align*}
The following sets will also be useful:
\begin{align*}
W_+^* &= W_+ \cap \{ x \in \R^4 \; : \; x_3 \neq 0\} \text{ and}
\\
W_-^* &= W_- \cap \{ x \in \R^4 \; : \; x_3 \neq 0\}.
\end{align*}
Now we consider the different life stages of orbits in $W^u(0)$.

\subsection*{Early life}
The next lemma gives sufficient control over the orbits in $W^u(0)$ early in their life when their dynamics are well approximated by the linearization of \eqref{FirstOrderODE} at the origin.

\begin{lemma}\label{EarlyLifeLem}
Suppose that $\Phi^0 : (-\infty, s_{\max})$ is a non-trivial orbit in $W^u(0)$ and $\sigma > 0$.
Then there exists an $s_0 \in (-\infty, s_{\max})$ and an $s$-translation of $Y$, denoted by $Y^0$, such that either:
\begin{enumerate}
\item
$\displaystyle
\left\{
\begin{aligned}
&\Phi^0_1(s_0) = Y^0_1(s_0),
\\
&\Phi^0(s_0) - Y^0(s_0) \in W_+^* \cup W_-^*, \text{ and}
\\
&|\Phi^0_3(s_0) - Y^0_3(s_0)| < \sigma;
\end{aligned}
\right.
$\vspace{1.0ex}%
\item
$\displaystyle
\left\{
\begin{aligned}
&-\Phi^0_1(s_0) = Y^0_1(s_0),
\\
&-\Phi^0(s_0) - Y^0(s_0) \in W_+^* \cup W_-^*, \text{ and}
\\
&|\Phi^0_3(s_0) + Y^0_3(s_0)| < \sigma;
\end{aligned}
\right.
$\vspace{1.0ex}%
\item
$\Phi^0(s_0) = Y^0(s_0)$; or
\vspace{1.0ex}%
\item
$\Phi^0(s_0) = -Y^0(s_0)$.
\qedhere
\end{enumerate}
\end{lemma}

\begin{pf}
Via the Stable Manifold theorem, see \cite[Section 2.7]{Perko91} for a proof, $W^u(0)$ is a smooth 2-manifold embedded in $\R^4$.

The linearization of \eqref{FirstOrderODE} at $\Phi = 0$ is
\begin{equation}\label{LinAtOriginEq}
\ps \Phi_{\lin}
=
\begin{pmatrix}
0 & 1 & 0 & 0 \\
0 & 0 & 1 & 0 \\
0 & 0 & 0 & 1 \\
-9 & 0 & 10 & 0
\end{pmatrix}
\Phi_{\lin}
=:
A \Phi_{\lin}.
\end{equation}
The eigenvalues of $A$ are $-3, -1, 1$, and $3$ with the corresponding eigenvectors
\begin{equation*}
\begin{pmatrix}
1\\
-3\\
9\\
-27
\end{pmatrix},
\begin{pmatrix}
1\\
-1\\
1\\
-1
\end{pmatrix},
\begin{pmatrix}
1\\
1\\
1\\
1
\end{pmatrix},
\mbox{ and }
\begin{pmatrix}
1\\
3\\
9\\
27
\end{pmatrix}.
\end{equation*}
Therefore, the tangent plane of $W^u(0)$ at the origin coincides with the linear subspace spanned by $(1,1,1,1)^T$ and $(1,3,9,27)^T$, hence $W^u(0)$ may be locally written as a graph over the $\Phi_1$-$\Phi_3$ plane:
\begin{equation}\label{LocalUnstableManifoldGraphFuncs}
\left\{
\begin{aligned}
\Phi_2(\Phi_1, \Phi_3) &= \frac{3}{4} \Phi_1 + \frac{1}{4} \Phi_3 + G_2(\Phi_1, \Phi_3),
\\
\Phi_4(\Phi_1, \Phi_3) &= -\frac{9}{4} \Phi_1 + \frac{13}{4} \Phi_3 + G_4(\Phi_1, \Phi_3),
\end{aligned}
\right.
\end{equation}
where
\begin{equation*}
\frac{\p(G_2, G_4)}{\p(\Phi_1, \Phi_3)}(0, 0) = 0.
\end{equation*}
Observe that, since $\Phi^0$ is not the trivial orbit, the dynamics of \eqref{FirstOrderODE} give, for any $s_0 \in (-\infty, s_{\max})$, an $s \in (-\infty, s_0]$ such that $\Phi^0_1(s) \neq 0$.
Therefore, we may find an $s_0$ sufficiently negative so that $\Phi^0_1(s_0) \neq 0$ and $|\Phi^0(s_0)|$ is as small as we wish.

Next, we assume that $\Phi^0_1(s_0) > 0$ and $s_0$ is sufficiently negative.
We take $Y^0$ to be an $s$-translation of the heteroclinic orbit $Y$ such that $Y^0_1(s_0) = \Phi^0_1(s_0)$.
Note that $Y^0$ is also an orbit in $W^u(0)$ and may be parameterized by $\Phi_1$.
If $\Phi^0_3(s_0) = Y^0_3(s_0)$ then $\Phi^0(s_0) = Y^0(s_0)$, since locally around the origin $W^u(0)$ is a graph over the $\Phi_1-\Phi_3$ plane.
This is Case (3) from the statement of this lemma.

From \eqref{LocalUnstableManifoldGraphFuncs}, we have
\begin{equation}\label{LocalUnstableManifoldGraphFuncsDerivatives}
\left\{
\begin{aligned}
\p_{\Phi_3} \Phi_2(\Phi_1, \Phi_3)
&=
\frac{1}{4} + o_{(\Phi_1, \Phi_3) \rightarrow 0}(1),
\\
\p_{\Phi_3} \Phi_4(\Phi_1, \Phi_3)
&=
\frac{13}{4} + o_{(\Phi_1, \Phi_3) \rightarrow 0}(1).
\end{aligned}
\right.
\end{equation}
If $\Phi^0_3(s_0) > Y^0_3(s_0)$ then $\Phi^0(s_0) - Y^0(s_0) \in W_+^*$.
On the other hand, if $\Phi^0_3(s_0) < Y^0_3(s_0)$ then $\Phi^0(s_0) - Y^0(s_0) \in W_-^*$.

Since we can choose $s_0$ sufficiently negative so that $|\Phi^0_3(s_0)|$ and $|Y^0(s_0)|$ are as small as we like, we can arrange for $|\Phi^0_3(s_0) - Y^0_3(s_0)| < \sigma$.
This is Case (1) from the statement of this lemma.

Recall that \eqref{FirstOrderODE} is invariant under the transformation $\Phi \mapsto -\Phi$.
Therefore, if $\Phi^0_1(s_0) < 0$ then we may argue the same as above but with $-\Phi^0$ instead of $\Phi^0$.
This leads to cases (2) and (4) from the statement of this lemma.
There are no more cases to consider.
\end{pf}

Due to symmetry, it suffices to only consider the cases (1) and (3) of Lemma \ref{EarlyLifeLem}.

\subsection*{Mid life}
Now that we have $\Phi^0 - Y^0 \in W_+ \cup W_-$, we can approximate \eqref{FirstOrderODE} by a system of simpler ordinary differential inequalities.

Let $\Phi^0 : (-\infty, s_{\max}) \rightarrow \R^4$, $Y^0$, and $s_0$ be the same as in Lemma \ref{EarlyLifeLem} with $\Phi^0 - Y^0 \in W_+ \cup W_-$.
We set $X(s) = \Phi^0(s) - Y^0(s)$ for $s \in [s_0, s_{\max})$.
Note that $\ps X_i = X_{i+1}$ for $i \in \{1,2,3\}$.

Before we prove our next result, we need an estimate.

\begin{lemma}\label{QGradLem}
Let
\begin{equation*}
f(y) = \frac{1}{2} \sin(2y) (3 \cos(2y) - 2)
\end{equation*}
and
\begin{equation*}
Q(x; f(y)) = \frac{2f(y) + 9 \sin(2x)}{14 + 6 \cos(2x)},
\end{equation*}
for $x, y \in \R$.
Then there exists a $c_0 \in (0,1)$ such that $\p_x Q(x;f(y)) \leq c_0$ for all $x, y \in \R$.
\end{lemma}

For the proof see the Appendix.

\begin{lemma}\label{LinIneqODE}
If $X \in W_+$ then $\ps X_4 \geq 4 (X_3 - c_0 X_1)$, and if $X \in W_-$ then $\ps X_4 \leq 4 (X_3 - c_0 X_1)$, where $c_0 \in (0,1)$ is taken from Lemma \ref{QGradLem}.
\end{lemma}

\begin{pf}
Recall that $Y_3 = \frac{1}{2} \sin(2 Y_1)$.
We have
\begin{align*}
\ps X_4
&=
F_1(Y_1^0 + X_1, Y_3^0 + X_3) - F_1(Y_1^0,Y_3^0)
\\
&\quad
+
(X_2 + Y_2^0)^2 F_2(Y_1^0 + X_1, Y_3^0 + X_3),
\end{align*}
since $F_2(Y_1^0, Y_3^0) = 0$.

If $X \in W_+$ then $F_2(Y_1^0 + X_1, Y_3^0 + X_3) \geq 0$, hence
\begin{equation*}
\ps X_4
\geq
F_1(Y_1^0 + X_1, Y_3^0 + X_3) - F_1(Y_1^0,Y_3^0).
\end{equation*}
On the other hand, if $X \in W_-$ then $F_2(Y_1^0 + X_1, Y_3^0 + X_3) \leq 0$, hence
\begin{equation*}
\ps X_4
\leq
F_1(Y_1^0 + X_1, Y_3^0 + X_3) - F_1(Y_1^0,Y_3^0).
\end{equation*}
Next, we study $F_1(Y_1^0 + X_1, Y_3^0 + X_3) - F_1(Y_1^0, Y_3^0)$.
Note that
\begin{equation*}
F_1(Y_1^0, Y_3^0) = F_1\left(Y_1^0, \frac{1}{2} \sin(2 Y_1^0)\right) = F_1(Y_1^0).
\end{equation*}
We are interested in the curve in the $\Phi_1$-$\Phi_3$ plane such that $F_1(\Phi_1, \Phi_3) = F_1(Y_1^0)$ for given values of $Y_1^0 \in (0, \pi)$.

This curve can be written as a graph over $\Phi_1$, namely:
\begin{equation*}
Q(\Phi_1; F_1(Y_1^0)) = \frac{2F_1(Y_1^0) + 9 \sin(2 \Phi_1)}{14 + 6 \cos(2 \Phi_1)}.
\end{equation*}
Therefore,
\begin{align*}
&\quad
F_1(Y_1^0 + X_1, Y_3^0 + X_3) - F_1(Y_1^0, Y_3^0)
\\
&=
F_1(Y_1^0 + X_1, Y_3^0 + X_3) - F_1(Y_1^0 + X_1, Q(Y_1^0 + X_1;F_1(Y_1^0))).
\end{align*}
Lemma \ref{QGradLem} yields $\p_{\Phi_1} Q(\Phi_1; F_1(Y_1)) \leq c_0$. Therefore, if $X \in W_+$ then
\begin{equation*}
Q(Y_1^0 + X_1; F_1(Y_1^0)) \leq Y_3^0 + c_0 X_1,
\end{equation*}
and if $X \in W_-$ then
\begin{equation*}
Q(Y_1^0 + X_1; F_1(Y_1^0)) \geq Y_3^0 + c_0 X_1.
\end{equation*}
Therefore, for $X \in W_+$:
\begin{align*}
\ps X_4
&\geq
F_1(Y_1^0 + X_1, Y_3^0 + X_3) - F_1(Y_1^0 + X_1, Y_3^0 + c_0 X_1)
\\
&=
\int_{Y_3^0 + c_0 X_1}^{Y_3^0 + X_3} \p_{\Phi_3} F_1(Y_1^0 + X_1, \Phi_3) \; d \Phi_3
\\
&\geq
4 (X_3 - c_0 X_1),
\end{align*}
and for $X \in W_-$:
\begin{align*}
\ps X_4
&\leq
F_1(Y_1^0 + X_1, Y_3^0 + X_3) - F_1(Y_1^0 + X_1, Y_3^0 + c_0 X_1)
\\
&=
\int_{Y_3^0 + X_3}^{Y_3^0 + c_0 X_1} -\p_{\Phi_3} F_1(Y_1^0 + X_1, \Phi_3) \; d \Phi_3
\\
&\leq
4 (X_3 - c_0 X_1).
\qedhere
\end{align*}
\end{pf}

The following lemma tells us that $W_+^*$ and $W_-^*$ are positive invariant sets for $X$.

\begin{lemma}\label{InvSetLem}
If $X(s_0) \in W_+^*$ (resp. $X(s_0) \in W_-^*$) then $X(s) \in W_+^*$ (resp. $X(s) \in W_-^*$) for $s \geq s_0$ while $X$ exists.
\end{lemma}

\begin{pf}
This follows easily from Taylor's theorem, \eqref{FirstOrderODE}, and Lemma \ref{LinIneqODE}.
\end{pf}

We now show that the non-trivial orbits in $W^u(0)$ that are not $s$-translations of $\pm Y$ must exit the region $|\Phi_3| < 1$ in finite time.

\begin{lemma}\label{ExitsStripLem}
Suppose that:
\begin{enumerate}[(i)]
\item
$\Phi^0 : [s_0, s_{\max}) \rightarrow \R^4$ solves \eqref{FirstOrderODE};
\item
$\Phi^0(s_0) - Y^0(s_0) \in W_+^*$ (resp. $\Phi^0(s_0) - Y^0(s_0) \in W_-^*$), where $Y^0$ is an $s$-translation of $Y$; and
\item
$|\Phi^0_3(s_0)| < 1$.
\end{enumerate}
Then there exists an $s_1 \in (s_0,  s_{\max})$ such that $\Phi^0_3(s_1) = 1$ (resp. $\Phi^0_3(s_1) = -1$).
\end{lemma}

\begin{remark}
Observe that $s_1 < s_{\max}$, because while $|\Phi^0_3| \leq 1$ we can control the growth of $|\Phi^0|$.
\end{remark}

\begin{pf}
We set $X = \Phi^0 - Y^0$.
Firstly, assume that $X(s_0) \in W_+^*$.
Lemmas \ref{LinIneqODE} and \ref{InvSetLem} give $\ps^2 X_3 \geq 2 X_3(s_0)$, hence
\begin{equation*}
\Phi^0_3(s_0+s)
\geq
Y^0_3(s_0+s) + X_3(s_0) s^2
\geq
-\frac{1}{2} + X_3(s_0) s^2.
\end{equation*}
Therefore, there exists such an $s_1 \in \left(s_0, s_0 + \left(\frac{3}{2X_3(s_0)}\right)^{\frac{1}{2}}\right]$.

The argument for $X(s_0) \in W_-^*$ is exactly the same.
\end{pf}

\subsection*{Late life}
The next lemma shows that once an orbit $\Phi^0 : (-\infty,s_{\max})$ in $W^u(0)$ has exited the region $|\Phi_3| < 1$ then it blows up in finite time, and in the process $|\Phi^0_1|$ diverges to infinity.

\begin{lemma}\label{BlowupLem}
Suppose that $\Phi_0 \in \R^4$ such that $\Phi_{0;3} \geq 1$ and $\Phi_{0;4} \geq 0$.
Let $\Phi^0 : [0, s_{\max}) \rightarrow \R^4$ be the solution to \eqref{FirstOrderODE} such that $\Phi^0(0) = \Phi_0$.
Then $\Phi^0$ blows up in finite time, that is, $s_{\max} < \infty$ and
\begin{equation*}
\lim_{s \nearrow s_{\max}} |\Phi^0(s)| = \infty.
\end{equation*}
Moreover, $\Phi^0_1(s) \rightarrow \infty$ as $s \nearrow s_{\max}$.
\end{lemma}

\begin{pf}
Let $\sigma > 0$ be arbitrary.
First we consider the case where $\Phi^0_2(0) \geq \sigma$.

Observe that the set
\begin{equation}\label{BlowupEq1}
S_1 = \{ \Phi \in \R^4 \; : \; \Phi_2 \geq \sigma, \Phi_3 \geq 1, \Phi_4 \geq 0\}
\end{equation}
is positive invariant under the flow described by \eqref{FirstOrderODE}.
For $\Phi^0 \in S_1$, we have
\begin{equation}\label{ApproxODE}
\ps \Phi^0_4 = F(\Phi^0) \left(\Phi_2^0\right)^2 \Phi_3^0,
\end{equation}
where $F(\Phi) \in [c_0, c_1] \subset (0,\infty)$, $c_0 = c_0(\sigma)$, and $c_1 = c_1(\sigma)$.

Observe that in \eqref{ApproxODE} $\Phi^0_1$ does not play a significant role.
Now we consider rescaled versions of $\Phi^0_2$ and $\Phi^0_3$:
\begin{equation*}
z_1 = \frac{\Phi^0_2}{(\Phi_4^0)^{\frac{1}{3}}} \mbox{ and }
z_2 = \frac{\Phi^0_3}{(\Phi_4^0)^{\frac{2}{3}}}.
\end{equation*}
We differentiate:
\begin{equation}\label{eq:zSystem}
\left\{
\begin{aligned}
\ps z_1
&=
\left(\Phi_4^0\right)^{\frac{1}{3}}
z_2
\left(
1
-
\frac{F(\Phi^0)}{3} z_1^3
\right) \mbox{ and}
\\
\ps z_2
&=
\left(\Phi_4^0\right)^{\frac{1}{3}}
\left(
1
-
\frac{2 F(\Phi^0)}{3}
 z_1^2 z_2^2
\right).
\end{aligned}
\right.
\end{equation}
Now \eqref{ApproxODE} becomes
\begin{equation*}
\ps \Phi^0_4 = F(\Phi^0) z_1^2 z_2  \left(\Phi_4^0\right)^{\frac{4}{3}}.
\end{equation*}
Problems arise with $z_1(0)$ and $z_2(0)$, if $\Phi_4^0(0) = 0$.
In this case we would like to examine $z_1(s)$ and $z_2(s)$ for $0 < s \ll 1$.
We have $\Phi^0(0) \in S_1$, hence $\Phi^0_4(s) > 0$ for all $s \in (0,s_{\max})$.
Therefore, $z_1$ and $z_2$ are well defined for $s \in (0,s_{\max})$, and $z_1(s), z_2(s) \rightarrow \infty$ as $s \searrow 0$.
On the other hand, if $\Phi^0_4(0) > 0$ then $z_1(0), z_2(0) > 0$.

Therefore, there exists an $\t{s} \in [0,s_{\max})$ such that $\Phi^0(\t{s}) \in S_1$, $\Phi^0_4(\t{s}) > 0$, and $z_1(\t{s}), z_2(\t{s}) > 0$.
We $s$-translate so that $\t{s} = 0$, and set
\begin{equation*}
Z = [z_{1;a}, z_{1;b}] \times [z_{2;a}, z_{2;b}] \subset (0,\infty) \times (0,\infty),
\end{equation*}
where
\begin{equation*}
z_{1;a}
=
\min\left\{z_1(0), \frac{1}{2} \left(\frac{3}{c_1}\right)^{\frac{1}{3}}\right\}, \;
z_{1;b}
=
\max\left\{z_1(0),2 \left(\frac{3}{c_0}\right)^{\frac{1}{3}}\right\},
\end{equation*}
and
\begin{equation*}
z_{2;a}
=
\min\left\{z_2(0), \frac{1}{2}\left(\frac{3}{2c_1}\right)^{\frac{1}{2}} z_{1;b}^{-1}\right\}, \;
z_{2;b}
=
\max\left\{z_2(0),2 \left(\frac{3}{2c_0}\right)^{\frac{1}{2}} z_{1;a}^{-1}\right\}.
\end{equation*}
Observe that $Z$ is a positive invariant set for \eqref{eq:zSystem}, and $(z_1(0), z_2(0)) \in Z$.
Therefore, for $s \in [0,s_{\max})$, we have
\begin{align*}
c_0 z_{1;a}^2 z_{2;a} \left(\Phi^0_4\right)^{\frac{4}{3}}
&\leq
\ps \Phi^0_4
\leq
c_1 z_{1;b}^2 z_{2;b} \left(\Phi^0_4\right)^{\frac{4}{3}},
\\
z_{1;a} \left(\Phi^0_4\right)^{\frac{1}{3}}
&\leq
\Phi^0_2
\leq
z_{1;b} \left(\Phi^0_4\right)^{\frac{1}{3}}, \mbox{ and}
\\
z_{2;a} \left(\Phi^0_4\right)^{\frac{2}{3}}
&\leq
\Phi^0_3
\leq
z_{2;b} \left(\Phi^0_4\right)^{\frac{2}{3}}.
\end{align*}
Therefore, $\Phi^0_4$ controls $|\Phi^0|$.
We have $\ps \Phi^0_4 \geq C \left(\Phi^0_4\right)^{\frac{4}{3}}$ and $\Phi^0_4(0) > 0$, hence $\Phi^0_4$ diverges to infinity in finite time, that is, $\Phi^0$ blows up in finite time.

Next, we turn our attention to showing that $\Phi^0_1 \rightarrow \infty$ as $s \nearrow s_{\max}$.
We let $i_0 \in \N$ be such that $2^{i_0} > \Phi^0_4(0)$.
For $i \in \N_0$, we let $s_i$ be defined via $\Phi^0_4(s_i) = 2^{i_0+i}$.
Since $\Phi_0^4$ is monotone increasing and diverges to infinity, these times are well-defined.
Because $\ps \Phi^0_4 \leq C (\Phi^0_4)^{\frac{4}{3}}$, we have that $s_{i+1} - s_i \geq C 2^{-\frac{1}{3}(i_0+i)}$.
For $s \in [s_i, s_{i+1}]$, we have $\ps \Phi_1^0(s) \geq  C 2^{\frac{1}{3}(i_0+i)}$.
Therefore, $\Phi_1^0(s_{i+1}) \geq \Phi_1^0(s_i) + C$ which implies $\Phi_1^0(s) \rightarrow \infty$ as $s \nearrow s_{\max}$, since $\Phi^0_1$ is monotone increasing.

Finally, we consider the case in which $\Phi_2(0) \leq 0$.
Observe that the set
\begin{equation}\label{BlowupEq2}
S_2 = \{ \Phi \in \R^4 \; : \; \Phi_3 \geq 1, \Phi_4 \geq 0\}
\end{equation}
is positive invariant under the flow described by \eqref{FirstOrderODE}.
Therefore, while $\Phi_2^0 \leq 0$ we have $|\Phi_2^0| \leq -\Phi_2^0(0)$ and
\begin{equation*}
|\ps \Phi_4^0| \leq C \Phi_3^0 (1 + |\Phi_2^0(0)|^2).
\end{equation*}
Therefore, there exists an $\t{s} \in (0, s_{\max})$ such that $\Phi^0_2(\t{s}) > 0$.
Now by autonomy we may $s$-translate, and then apply the previous argument to the new initial data $\Phi^0(\t{s})$.
\end{pf}

%%%%%%%%%%%%%%%%%%%%%%%%%%%%%%%%%%%%%%%%%%%%%%%%%%%%%%%%%%%%
\subsection*{Proof of Theorem \ref{MainUniqueThm}}
%%%%%%%%%%%%%%%%%%%%%%%%%%%%%%%%%%%%%%%%%%%%%%%%%%%%%%%%%%%%
If $\Phi_0 = Y^0$ or $-\Phi^0 = Y^0$, where $Y^0$ is an $s$-translation of $Y$, then we are done.
Otherwise, Lemma \ref{EarlyLifeLem} tells us that at some $s_0 \in (-\infty, s_{\max})$ we have $\Phi^0(s_0) - Y_0(s_0) \in W^*_+ \cup W^*_-$ or $-\Phi^0(s_0) - Y^0(s_0) \in W^*_+ \cup W^*_-$.
Since \eqref{FirstOrderODE} is invariant under the transformation $\Phi \mapsto -\Phi$, it suffices to only consider the case where $\Phi^0(s_0) - Y^0(s_0) \in W^*_+ \cup W^*_-$.
Now Lemma \ref{ExitsStripLem} shows that there exists an $s_1 \in (-\infty, s_{\max})$ such that $|\Phi^0_3(s_1)| = 1$.
Again due to the invariance of \eqref{FirstOrderODE} under $\Phi \mapsto -\Phi$, we may assume that $\Phi^0_3(s_1) = 1$.
Finally, Lemma \ref{BlowupLem} shows that $\Phi^0$ must blowup in finite time.
%%% Local Variables:
%%% mode: latex
%%% TeX-master: "main"
%%% End:

\section{Further properties of the unstable manifold}\label{UnstableManifoldSect}
\makeatletter{}In this section we continue our study of $W^u(0)$.
This leads to the proofs of theorems \ref{MainNonExistence} and \ref{MainArbitraryWrapAround}.

Our first result concerns non-trivial orbits $\Phi^0$ in $W^u(0)$ which are not $s$-translations of $\pm Y$.
We know that these orbits must exit the region $|\Phi_3| < 1$ in finite time.
The next lemma tells us that $(\Phi^0_1, \Phi^0_2)$ stays within a bounded region of $\R^2$ up until and including this exit time.

\begin{lemma}\label{BoundedSpeedAtExitLem}
Let $\Phi^0 : (-\infty, s_{\max}) \rightarrow \R^4$ be an orbit in $W^u(0)$ such that $|\Phi^0_3(\t{s})| = 1$ for some $\t{s} \in (-\infty, s_{\max})$.
Moreover, let $\t{s}$ be the first such time in which $|\Phi^0_3(\t{s})| = 1$.
Then $|(\Phi^0_1(s), \Phi^0_2(s))| \leq C$ for all $s \in (-\infty, \t{s}]$.
\end{lemma}

\begin{pf}
Since $\Phi^0$ is non-trivial and not an $s$-translation of $\pm Y$, we are in either Case (1) or Case (2) of Lemma \ref{EarlyLifeLem}.

For now we assume that we are in Case (1).
This means that we have an $s_0 \in (-\infty,s_{\max})$, a $Y^0$ which is an $s$-translation of $Y$, and an $X = \Phi^0 - Y^0$ such that $X(s_0) \in W^*_+ \cup W^*_-$ and $X_3(s_0) \neq 0$ is as small as we like, in particular $|X_3(s_0)| \leq \frac{1}{4}$.
Due to autonomy we may assume that $s_0 = 0$.
Observe that \eqref{LocalUnstableManifoldGraphFuncsDerivatives} gives
\begin{equation}\label{BoundedSpeedAtExitEq1}
|X_2(0)| \leq C |X_3(0)|
\mbox{ and }
|X_4(0)| \leq C |X_3(0)|.
\end{equation}
On intervals on which we have uniform control of $X_3$ we also have uniform control of $\Phi_3$.
Therefore, blowup may not happen on such intervals.
In what follows we study $X_3$ for $|X_3| \leq 2$.

We let $i_0$ be the largest integer such that $2^{-i_0+1} > |X_3(0)|$.
We set $s_0 = 0$ and $|X_3(s_i)| = 2^{-i_0 + i}$ for $i \in \N$.
For $s \in [s_i, s_{i+1}]$, we have $\ps |X_4| \geq 2^{-i_0 + i}$, hence $|X_3(s)| \geq |X_3(s_i)| + 2^{-i_0+i-1} (s_{i+1}-s_i)^2$ which implies that $s_{i+1} - s_i \leq C$.

For $s \in [s_i, s_{i+1}]$, we have $\ps |X_2| \leq 2^{-i_0 + i +1}$, hence
\begin{equation*}
|X_2(s)|
\leq
|X_2(0)|
+
C \sum_{j=0}^{i-1} 2^{-i_0+j+1}
\leq
C
\left(|X_3(0)| + 2^{-i_0 + i + 1}\right)
\mbox{ for } s \in [s_0, s_i],
\end{equation*}
where we have used \eqref{BoundedSpeedAtExitEq1}.
Observe that $\t{s} \in [0, s_{i_0+1}]$, hence $|X_2(s)| \leq C$ for all $s \in (-\infty, \t{s}]$.
This implies that, for all $s \in (-\infty, \t{s}]$, $|\Phi^0_2(s)| \leq C$.
Moreover, since $X(s) \in W_+ \cup W_-$, we have $|\Phi^0_1(s)| \leq C$.

Finally, we look at what happens if we are in Case (2).
Because of the invariance of \eqref{FirstOrderODE} under the transformation $\Phi \mapsto -\Phi$, we may use the above argument on $-\Phi^0$ yielding the same conclusion of $|(\Phi^0_1(s), \Phi^0_2(s))| \leq C$ for all $s \in (-\infty, \t{s}]$.
\end{pf}

Using this result, we show that, for an equivariant map from $\ol{B^4(0, 1)}$ into $S^4$, if the normal derivative at the boundary vanishes then there is a limit on the number of times an equivariant biharmonic map from $\ol{B^4(0,1)}$ into $S^4$ satisfying the same boundary condition can wind around $S^4$.

\begin{lemma}\label{NonExistenceLem}
Let $\Phi \in W^u(0)$.
Then there exists a $K > 0$ such that if $\Phi_2 = 0$ then $|\Phi_1| \leq K$.
\end{lemma}

\begin{pf}
Let $\Phi^0 : (-\infty,s_{\max}) \rightarrow \R^4$ be an orbit in $W^u(0)$.
If $\Phi^0$ is trivial or an $s$-translation of $\pm Y$ then $|\Phi^0_1(s)| < \pi$ for all $s \in (-\infty, s_{\max})$, hence if we take $K > \pi$ then these cases cause us no issues.

Otherwise, we find ourselves in Case (1) or Case (2) of Lemma \ref{EarlyLifeLem}.
In these cases Lemma \ref{ExitsStripLem} gives a time $\t{s} \in (-\infty, s_{\max})$ such that $|\Phi^0_3(\t{s})| = 1$.
We may assume that $\t{s}$ is the first such time.
Lemma \ref{BoundedSpeedAtExitLem} shows that $|(\Phi^0_1(s), \Phi^0_2(s))| \leq C$ for all $s \in (-\infty,\t{s}]$.

Let us assume that $\Phi^0_3(\t{s}) = 1$.
If $\Phi^0_2(\t{s}) > 0$ then it will be positive for all $s \in [\t{s},s_{\max})$, since $S_1$ from \eqref{BlowupEq1} is a positive invariant set.
Therefore, if $\Phi^0_2(s) = 0$ then $s < \t{s}$ and $|\Phi^0_1(s)| \leq C$.
On the other hand, if $\Phi^0_2(\t{s}) \leq 0$ then $\ps \Phi^0_2(s) \geq 1$ for all $s \in [\t{s},s_{\max})$, since $S_2$ from \eqref{BlowupEq2} is a positive invariant set.
Therefore, while $\Phi^0_2(s) \leq 0$ we have $|\Phi^0_2(s)| \leq |\Phi^0_2(\t{s})|$ for $s \in [\t{s},s_{\max})$.
If $\Phi^0_2(s) < 0$ for all $s \in [\t{s}, s_{\max})$ then there is nothing more to consider.
On the other hand, there is a unique $s_0 \in [\t{s}, s_{\max})$ such that $\Phi^0_2(s_0) = 0$.
Observe that $s_0 - \t{s} \leq |\Phi^0_2(\t{s})| \leq C$, hence
\begin{equation}\label{NonExistenceEq1}
|\Phi^0_1(s_0)| \leq |\Phi^0_1(\t{s})| + (s_0 - \t{s}) |\Phi^0_2(\t{s})| \leq C.
\end{equation}
Therefore, if we take $K > 0$ sufficiently large, these cases also do not cause us any problems.

Finally, we consider the case in which $\Phi^0_3(\t{s}) = -1$.
Due to the invariance of \eqref{FirstOrderODE} under the transformation $\Phi \mapsto -\Phi$, we may apply the above argument to $-\Phi^0$.
\end{pf}

Theorem \ref{MainNonExistence} is a corollary of this.

Next, we show the existence of smooth equivariant biharmonic maps from $B^4(0,1)$ into $S^4$ that can wind around $S^4$ as many times as we wish.
Before we do this we need some preparatory lemmas.
Our arguments are influenced by the ideas in \cite{GastelZorn12}.
Recall that given an orbit $\Phi^0 : (-\infty,0] \rightarrow \R^4$ in $W^u(0)$, $\psi(r) = \Phi^0_1(\log r)$ solves \eqref{BMHFSymRedEq} (with $\pt \psi = 0$) on $(0,1]$.
Note that $\psi(0) = 0$.
We first wish to verify that given such a $\psi$, $u = \Upsilon(\psi)$ is weakly biharmonic.

The next lemma obtains estimates on the derivatives of our solutions $\psi$ and the corresponding equivariant maps.

\begin{lemma}\label{OriginGrowthLem}
Let $\psi \in C([0,1];\R) \cap C^{\infty}((0,1];\R)$, with $\psi(0) = 0$, be a solution to \eqref{BMHFSymRedEq} (with $\pt \psi = 0$) and
\begin{equation*}
u = \Upsilon(\psi) \in C(\ol{B^4(0,1)}; S^4) \cap C^{\infty}(\ol{B^4(0,1)} \setminus \{0\};S^4).
\end{equation*}
Then, for $r > 0$,
\begin{equation*}
|\psi(r)|
\leq
C r, \;
|\pr \psi(r)|
\leq
C, \mbox{ and }
|\pr^2 \psi(r)|
\leq
Cr.
\end{equation*}
Furthermore, $Du \in L^{\infty}(B^4(0,1))$ and $|D^2 u(x)| \leq C |x|^{-1}$ for $x \in \ol{B^4(0,1)} \setminus \{0\}$.
In the above inequalities $C = C(\psi)$.
\end{lemma}

\begin{pf}
We set $\phi(s) = \psi(e^s)$ and $\Phi^0_i = \ps^{i-1} \phi$ for $i \in \{1,2,3,4\}$.
Recall that $\Phi^0$ solves \eqref{FirstOrderODE}.

We rewrite \eqref{FirstOrderODE} as
\begin{equation}\label{BdryCondUnstableEq2}
\ps \Phi^0 = A \Phi^0 + G(\Phi^0),
\end{equation}
where $A$ is the same as in \eqref{LinAtOriginEq} and $|G(\Phi^0)| \leq C |\Phi^0|^3$ for sufficiently small $\Phi^0$.
We set
\begin{equation}\label{ChangeOfVariables}
V^0(s)
=
\begin{pmatrix}
1 & 1  & 1   & 1\\
1 & 3  & -1  & -3\\
1 & 9  & 1   & 9\\
1 & 27 & -1  & -27
\end{pmatrix}^{-1}
\Phi^0(-s)
=:
P^{-1} \Phi^0(-s).
\end{equation}
We substitute this into \eqref{BdryCondUnstableEq2}:
\begin{equation}\label{TransformedFirstOrderODE}
\ps V^0
=
\begin{pmatrix}
-1 & 0 & 0 & 0 \\
0 & -3 & 0 & 0 \\
0 & 0 & 1 & 0 \\
0 & 0 & 0 & 3
\end{pmatrix}
V^0
-
P^{-1} G(P V^0)
=:
D V^0
+
\t{G}(V^0),
\end{equation}
where $|\t{G}(V)| \leq C |V|^3$ for sufficiently small $V$.
Since we reversed $s$ in \eqref{ChangeOfVariables}, we are now interested in the stable manifold at the origin of \eqref{TransformedFirstOrderODE}.
This manifold is tangent to the $V_1$ -- $V_2$ plane at the origin, and can be locally written as a graph over this plane with $V_3 = V_3(V_1, V_2)$ and $V_4 = V_4(V_1, V_3)$ such that
\begin{equation*}
\frac{\p (V_3, V_4)}{\p (V_1, V_2)}(0,0) = 0.
\end{equation*}
Our first aim is to show that
\begin{equation}\label{BdryCondUnstableEq5}
|(V_3(V_1, V_2), V_4(V_1, V_2))| \leq C |(V_1, V_2)|^3,
\end{equation}
for sufficiently small $|(V_1, V_2)|$.

We let $\eps > 0$, and $(V_1, V_2) \in \R^2$ such that $|(V_1, V_2)| \leq \eps$.
We setup the iteration:
\begin{equation*}
\ps V^1 = D V^1 \mbox{ with } V^1(0) = (V_1, V_2, 0, 0),
\end{equation*}
and for $i \in \N$:
\begin{equation*}
\left\{
\begin{aligned}
\ps V^{i+1} &= D V^{i+1} + \t{G}(V^i),
\\
V^{i+1}(0) &= (V_1, V_2, V_{3;i+1}, V_{4;i+1}),
\end{aligned}
\right.
\end{equation*}
where
\begin{equation}\label{BdryCondUnstableEq4}
\left\{
\begin{aligned}
V_{3;i+1}
&=
-\int_0^{\infty} e^{-s} \t{G}_3(V^i(s)) \ds, \mbox{ and}
\\
V_{4;i+1}
&=
-\int_0^{\infty} e^{-3s} \t{G}_4(V^i(s)) \ds.
\end{aligned}
\right.
\end{equation}
This iteration is used in \cite[Section 2.7]{Perko91} as part of the proof of the Stable Manifold theorem.
In this proof it is shown that, for sufficiently small $\eps$,
\begin{equation*}
V_3(V_1, V_2) = \lim_{i \rightarrow \infty} V_{3;i} \text{ and } V_4(V_1, V_2) = \lim_{i \rightarrow \infty} V_{4;i}.
\end{equation*}
Furthermore, it is shown that
\begin{equation}\label{BdryCondUnstableEq3}
|V^i(s)| \leq C |(V_1, V_2)| e^{-\alpha s},
\end{equation}
for any $-\alpha > -1$, as long as $\eps$ is sufficiently small.
We substitute \eqref{BdryCondUnstableEq3} into \eqref{BdryCondUnstableEq4}, and take the limit $i \rightarrow \infty$ to obtain \eqref{BdryCondUnstableEq5}.

If $\Phi^0$ is an orbit in $W^u(0)$ then $V^0$ is an orbit in the stable manifold of the origin of \eqref{TransformedFirstOrderODE}.
For sufficiently large $s$, \eqref{TransformedFirstOrderODE} and \eqref{BdryCondUnstableEq5} give
\begin{equation}\label{OriginGrowthVEq}
|V^0_1(s)| \leq C e^{-s},\,
|V^0_2(s)| \leq C e^{-3s},\,
|V^0_3(s)| \leq C e^{-3s}, \text{ and }
|V^0_4(s)| \leq C e^{-3s}.
\end{equation}
For $r > 0$, we have
\begin{equation*}
\pr \psi(r)
=
\frac{\Phi^0_2(\log r)}{r}, \mbox{ and }
\pr^2 \psi(r)
=
\frac{\Phi^0_3(\log r) - \Phi^0_2(\log r)}{r^2}.
\end{equation*}
Using this, along with \eqref{ChangeOfVariables} and \eqref{OriginGrowthVEq}, gives us
\begin{equation}\label{OriginGrowthEqAll}
|\psi(r)| \leq C(\psi) r,
\;
|\pr \psi(r)| \leq C(\psi),
\mbox{ and }
|\pr^2 \psi(r)|
\leq
C(\psi) r.
\end{equation}
Next, we turn our attention towards the estimates on $u$.
First we focus on $|Du|$.
Recalling the notation from \eqref{BraceNotationEq}, for $x \neq 0$, we have
\begin{equation*}
D (\{f_0, f_1\}_{\h{x}})
:
D (\{g_0, g_1\}_{\h{x}})
=
\pr f_0 \; \pr g_0 + \pr f_1 \; \pr g_1
+
\frac{d-1}{r^2} f_0 \, g_0.
\end{equation*}
This, \eqref{OriginGrowthEqAll}, and $u = \Upsilon(\psi)$ yield $Du \in L^{\infty}(B^4(0,1))$.
For $d = 4$, we have
\begin{equation*}
\begin{aligned}
|D^2 \{f_0, f_1\}_{\h{x}}|^2
&=
\frac{1}{r^4}
\Big(
9 f_0^2-18 r f_0 \pr f_0 +r^2 \Big(9 (\pr f_0)^2+3 (\pr f_1)^2
\\
&\quad
+r^2 \left((\pr^2 f_0)^2+(\pr^2 f_1)^2\right)\Big)
\Big).
\end{aligned}
\end{equation*}
This, \eqref{OriginGrowthEqAll}, and $u = \Upsilon(\psi)$ yield $|D^2 u| \leq C(\psi) |x|^{-1}$.
\end{pf}

Now we wish to show that our solutions $\psi$ to \eqref{BMHFSymRedEq} (with $\pt \psi = 0$), with $\psi(0) = 0$, give rise to equivariant maps $u \in H^2$.

\begin{lemma}\label{InALem}
Let $\psi \in C([0,1];\R) \cap C^{\infty}((0,1];\R)$, with $\psi(0) = 0$, be a solution to \eqref{BMHFSymRedEq} (with $\pt \psi = 0$) and
\begin{equation*}
u = \Upsilon(\psi) \in C(\ol{B^4(0,1)}; S^4) \cap C^{\infty}(\ol{B^4(0,1)} \setminus \{0\}; S^4).
\end{equation*}
Then $u \in H^2(B^4(0,1);S^4)$.
\end{lemma}

\begin{pf}
Using the growth estimates from Lemma \ref{OriginGrowthLem}, it can be shown that $u$ has weak derivatives in $L^2$ up to order two which are equal a.e.\ to their respective classical derivatives.
\end{pf}

Next, we show that our solutions $\psi$ give rise to weakly biharmonic maps.
See \cite{GastelZorn12} for a different approach in a slightly different situation.

\begin{lemma}\label{FirstVariationVerificationLem}
Let $\psi \in C([0,1];\R) \cap C^{\infty}((0,1];\R)$, with $\psi(0) = 0$, be a solution to \eqref{BMHFSymRedEq} (with $\pt \psi = 0$), and
\begin{equation*}
u = \Upsilon(\psi) \in C(\ol{B^4(0,1)}; S^4) \cap C^{\infty}(\ol{B^4(0,1)} \setminus \{0\}; S^4).
\end{equation*}
Then $u$ is weakly biharmonic.
\end{lemma}

\begin{pf}
We let $\eta \in C^{\infty}_c(B^4(0,1); \R^5)$ be arbitrary.
We wish to show that
\begin{equation*}
\pt|_{t = 0} E_2(\Pi(u + t \eta)) = 0,
\end{equation*}
where $\Pi(x) = \frac{x}{|x|}$ is defined on $\R^5 \setminus \{0\}$.

From \cite[(2.1) and (2.2)]{Strzelecki03}, we have
\begin{equation*}
\pt|_{t = 0} E_2(\Pi(u + t \eta))
=
2 \int_{B^4(0,1)} \left( \Delta u \cdot \Delta \eta - \sum_{\gamma=1}^5 \Delta u^{\gamma} \Delta \left( u^{\gamma} u \cdot \eta\right) \right) \dx.
\end{equation*}
We let $\omega \in C^{\infty}_c(B^4(0,1); [0,1])$ be such that $\omega \equiv 1$ on $B^4\left(0,\frac{1}{2}\right)$.
For $R > 0$, we set $\omega_R(x) = \omega(x/R)$.
We have
\begin{equation*}
\begin{aligned}
\pt|_{t = 0} E_2(\Pi(u + t \eta))
&=
\pt|_{t = 0} E_2(\Pi(u + t (\omega_R \eta)))
\\
&\quad
+
\pt|_{t = 0} E_2(\Pi(u + t ((1-\omega_R)\eta))).
\end{aligned}
\end{equation*}
Lemma \ref{OriginGrowthLem}, gives us
\begin{equation*}
\pt|_{t = 0} E_2(\Pi(u + t (\omega_R \eta)))
=
o_{R \searrow 0}(1).
\end{equation*}
Next, we turn our attention towards $\pt|_{t = 0} E_2(\Pi(u + t ((1-\omega_R)\eta)))$.
Since the support of $(1-\omega_R)\eta$ is bounded away from the origin, and $u$ is smooth and satisfies the Euler-Lagrange equation \eqref{BiEulerLagrangeEq} away from the origin, we have $\pt|_{t = 0} E_2(\Pi(u + t ((1-\omega_R)\eta))) = 0$.
Therefore, $\pt|_{t = 0} E_2(\Pi(u + t \eta)) = o_{R \searrow 0}(1)$ which gives the desired result after taking the limit $R \searrow 0$.
\end{pf}

Finally, we prove Theorem \ref{MainArbitraryWrapAround}.

%%%%%%%%%%%%%%%%%%%%%%%%%%%%%%%%%%%%%%%%%%%%%%%%%%%%%%%%%%%%
\subsection*{Proof of Theorem \ref{MainArbitraryWrapAround}}
%%%%%%%%%%%%%%%%%%%%%%%%%%%%%%%%%%%%%%%%%%%%%%%%%%%%%%%%%%%%
Due to the invariance of \eqref{FirstOrderODE} under the transformation $\Phi \mapsto -\Phi$, it suffices to prove the result for $a \geq 0$.
The $a = 0$ case is taken care of by the trivial solution.

Therefore, we let $a > 0$ be arbitrary.
There exists a non-trivial orbit
\begin{equation*}
\Phi^0 : (-\infty, s_{\max}) \rightarrow \R^4
\end{equation*}
in $W^u(0)$ which is not an $s$-translation of $\pm  Y$.
Lemma \ref{ExitsStripLem} tells us that $\Phi^0$ must exit the region $|\Phi_3| < 1$ in finite time.
Due to the invariance of \eqref{FirstOrderODE} under the transformation $\Phi \mapsto -\Phi$, we may assume that there exists some $s_0 \in (-\infty, s_{\max})$ such that $\Phi^0_3(s_0) = 1$.
Lemma \ref{BlowupLem} tells us that $\Phi^0_1(s) \rightarrow \infty$ as $s \nearrow s_{\max}$.
Since $\Phi^0$ is an orbit in $W^u(0)$, we also know that $\Phi^0_1(s) \rightarrow 0$ as $s \rightarrow -\infty$.
Therefore, we may $s$-translate $\Phi^0$ so that $\Phi^0_1(0) = a$.
This corresponds to a solution of \eqref{FourthOrderODE} with $\phi(0) = a$, which after undoing the change of coordinates $r = e^s$, corresponds to a solution of \eqref{BMHFSymRedEq} (with $\pt \psi = 0$) such that $\psi(1) = a$.

Lemmas \ref{InALem} and \ref{FirstVariationVerificationLem}, tell us that
\begin{equation*}
u = \Upsilon(\psi) \in C^{\infty}(\ol{B^4(0,1)} \setminus \{0\}; S^4) \cap H^2(B^4(0,1); S^4),
\end{equation*}
is a weakly biharmonic map.
Standard higher interior regularity arguments, see for example \cite{CWY99}, yield smoothness of $u$ on all of $\ol{B^4(0,1)}$.
%%% Local Variables:
%%% mode: latex
%%% TeX-master: "main"
%%% End:

% setup appendix
\setcounter{section}{0}
\setcounter{equation}{0}
\renewcommand{\theequation}{\thesection.\arabic{equation}}
\appendix
\section*{Appendix}\label{TechnicalSect}
\renewcommand{\thesection}{A}
\makeatletter{}\begin{pf}[Proof of Lemma \ref{BdryCondUnstableLem}]
{\bf(2) $\Longrightarrow$ (1):}
This direction is trivial.

{\bf(1) $\Longrightarrow$ (2):}
We set $\t{\Phi}^0_i(s) = (-1)^{i+1} \Phi^0_i(-s)$ for $i \in \{1,2,3,4\}$.
Observe that $\Phi^0$ solving \eqref{FirstOrderODE} is equivalent to $\t{\Phi}^0$ solving \eqref{FirstOrderODE}.
Therefore, after relabeling $\t{\Phi}^0$ as $\Phi^0$ our statement is equivalent to showing that if $s_0 \in \R$, $\Phi^0 : [s_0, \infty) \rightarrow \R^4$ solves \eqref{FirstOrderODE}, and
\begin{equation*}
\lim_{s \rightarrow \infty} \Phi^0_1(s) = 0,
\end{equation*}
then
\begin{equation*}
\lim_{s \rightarrow \infty} \Phi^0(s) = 0.
\end{equation*}
It is easy to show that if $x \in C^2([s_0,\infty);\R)$, $x(s) \rightarrow 0$ as $s \rightarrow \infty$, and $|\ps^2 x(s)| \leq C$ for all $s \in [s_0, \infty)$ then $\ps x(s) \rightarrow 0$ as $s \rightarrow \infty$.
We use this fact, which we call (P1), repeatedly in what follows.

% Boundedness of \Phi^0_3
First observe that there cannot exist an $s_1 \in [s_0, \infty)$ such that $|\Phi^0_3(s)| \geq 1$ for all $s \in [s_1, \infty)$.
Indeed, if there were such an $s_1$ then eventually $\Phi^0_4(s)$ would be the same sign as $\Phi^0_3$ after which we could apply Lemma \ref{BlowupLem} and obtain a contradiction.

Therefore, if there is an $s_1 \in [s_0, \infty)$ such that $|\Phi^0_3(s_1)| \geq 1$ then there must be an $s_2 > s_1$ such that $|\Phi^0_3(s)| < 1$ for all $s \in [s_2, \infty)$, or else we could apply Lemma \ref{BlowupLem} and obtain a contradiction.
Therefore, $\Phi^0_3$ is bounded on $[s_0, \infty)$.

% \Phi^0_2
Now we proceed to show, one by one, that $\lim_{s\rightarrow \infty} \Phi^0_i(s) = 0$ for $i \in \{2,3,4\}$.
First we look at $\Phi^0_2$.
The fact that $\Phi^0_1 \rightarrow 0$ as $s \rightarrow \infty$, the boundedness of $\Phi^0_3$, and (P1) yield
\begin{equation*}
\lim_{s \rightarrow \infty} \Phi^0_2(s) = 0.
\end{equation*}
% \Phi^0_3
Next, we look at $\Phi^0_3$.
Hoping for a contradiction, we assume that $\Phi^0_3(s) \not\rightarrow 0$ as $s \rightarrow \infty$.
From (P1) we know that $\Phi^0_4$ is unbounded, that is, there exists a monotone increasing sequence $\{s_i\}_{i \in \N} \subset [s_0, \infty)$ diverging to infinity such that $|\Phi^0_4(s_i)| \rightarrow \infty$.
From \eqref{FirstOrderODE} and the fact that $|(\Phi^0_1(s), \Phi^0_2(s), \Phi^0_3(s))| \leq C$ on $[s_0,\infty)$, we have that $|\ps \Phi^0_4(s)| \leq C$ on $[s_0, \infty)$.
Therefore, $|\Phi^0_4(s)| \geq \frac{1}{2} |\Phi^0_4(s_i)|$ for $s \in \left[s_i, s_i+\frac{1}{2C} |\Phi^0_4(s_i)|\right]$.
Observe that over this interval $\Phi^0_4$ is non-vanishing.
Therefore, there exists an $s \in [s_0,\infty)$ such that $|\Phi^0_3(s)| \geq 1$ and $\Phi^0_3(s)$ has the same sign as $\Phi^0_4(s) \neq 0$.
Lemma \ref{BlowupLem} then yields a contradiction, hence
\begin{equation*}
\lim_{s \rightarrow \infty} \Phi^0_3(s) = 0.
\end{equation*}
% \Phi^0_4
Finally, we look at $\Phi^0_4$.
Since
\begin{equation*}
\lim_{s \rightarrow \infty} (\Phi^0_1(s), \Phi^0_2(s), \Phi^0_3(s)) = 0,
\end{equation*}
from \eqref{FirstOrderODE}, we have $\ps \Phi^0_4(s) \rightarrow 0$ as $s \rightarrow \infty$.
Now (P1) gives us
\begin{equation*}
\lim_{s \rightarrow \infty} \Phi^0_4(s) = 0.
\qedhere
\end{equation*}
\end{pf}

\begin{remark}
Observe that in the above proof we make use of Lemma \ref{BlowupLem}.
Our argument would be circular if Lemma \ref{BlowupLem} depended upon Lemma \ref{BdryCondUnstableLem}.
By closely examining the proof of Lemma \ref{BlowupLem}, it is clear that this is not the case.
\end{remark}

\begin{pf}[Proof of Lemma \ref{QGradLem}]
We prove this lemma for $c_0 = \frac{99}{100}$. It is elementary to compute
\begin{equation*}
\min_{x \in \R} f(x) = -\frac{1}{12} \sqrt{169+38 \sqrt{19}}.
\end{equation*}
Since $f$ is an odd function, we have
\begin{equation*}
|f(y)| \leq \frac{1}{12} \sqrt{169+38 \sqrt{19}} \leq 2 \mbox{ for all } y \in \R.
\end{equation*}
We differentiate:
\begin{equation*}
\partial_x Q(x;f(y))
=
\frac{3 (9+21 \cos(2x)+2 f(y) \sin(2x))}{(7+3 \cos(2x))^2}.
\end{equation*}
By periodicity, what we wish to prove is that $\px Q(x;f(y)) \leq c_0$ for all $x \in \left(-\frac{\pi}{2}, \frac{\pi}{2}\right]$ and $y \in \R$.

Firstly, $\px Q\left(\frac{\pi}{2};f(y)\right) < 0$ which means we may restrict our attention to $x \in \left(-\frac{\pi}{2}, \frac{\pi}{2}\right)$.
We use Weierstrass' substitution:
\begin{equation*}
\sin(2x) \mapsto \frac{2t}{1+t^2} \mbox{ and } \cos(2x) \mapsto \frac{1-t^2}{1+t^2} \mbox{ for } t \in \R.
\end{equation*}
This transforms the problem into showing that
\begin{equation*}
\frac{3 (15+2 (\t{f}-3 t) t) \left(1+t^2\right)}{2 \left(5+2 t^2\right)^2}
\leq c_0,
\end{equation*}
for all $t \in \R$ and $\t{f} \in [-2,2]$.
It suffices to show
\begin{equation}\label{QGradEq1}
-45+50 c_0+(-27+40 c_0) t^2+(18+8 c_0) t^4-12 (t+ t^3) \geq 0,
\end{equation}
for all $t \in \R$.
Next we prove this.

We substitute $c_0 = \frac{99}{100}$ into \eqref{QGradEq1} and let $p$ be the polynomial on the left hand side of the resulting expression, that is,
\begin{equation*}
p(t)
=
\frac{648}{25} t^4-12 t^3+\frac{63}{5} t^2-12 t+\frac{9}{2}.
\end{equation*}
We calculate:
\begin{equation*}
p'\left(\frac{2}{5}\right) < 0,
\;
p'\left(\frac{43}{100}\right) > 0,
\mbox{ and }
p''(t) > 0,
\end{equation*}
for $t \in \R$.
Therefore, $p$ is convex with its unique global minimum occurring somewhere in $\left[\frac{2}{5}, \frac{43}{100}\right]$.
We use this to estimate:
\begin{equation*}
\min_{t \in \R} p(t)
\geq
\frac{648}{25} \left(\frac{2}{5}\right)^4-12 \left(\frac{43}{100}\right)^3+\frac{63}{5} \left(\frac{2}{5}\right)^2-12 \left(\frac{43}{100}\right)+\frac{9}{2}
>
0.
\end{equation*}
This is what we wished to show.
\end{pf}
%%% Local Variables:
%%% mode: latex
%%% TeX-master: "main"
%%% End:

\bibliographystyle{plain}
\bibliography{master}

\begin{thebibliography}{10}

\bibitem{Angelsberg06}
G.~Angelsberg.
\newblock A monotonicity formula for stationary biharmonic maps.
\newblock {\em Math. Z.}, 252(2):287--293, 2006.

\bibitem{ChangDing91}
K.-C. Chang and W.Y. Ding.
\newblock A result on the global existence for heat flows of harmonic maps from
  {$D^2$} into {$S^2$}.
\newblock In {\em Nematics ({O}rsay, 1990)}, volume 332 of {\em NATO Adv. Sci.
  Inst. Ser. C Math. Phys. Sci.}, pages 37--47. Kluwer Acad. Publ., Dordrecht,
  1991.

\bibitem{ChangDingYe1992}
K.-C. Chang, W.Y. Ding, and R.~Ye.
\newblock Finite-time blow-up of the heat flow of harmonic maps from surfaces.
\newblock {\em J. Differential Geom.}, 36(2):507--515, 1992.

\bibitem{ChangGurskyYang99}
S.-Y.~A. Chang, M.J. Gursky, and P.C. Yang.
\newblock Regularity of a fourth order nonlinear {PDE} with critical exponent.
\newblock {\em Amer. J. Math.}, 121(2):215--257, 1999.

\bibitem{CWY99}
S.-Y.~A. Chang, L.~Wang, and P.C. Yang.
\newblock A regularity theory of biharmonic maps.
\newblock {\em Comm. Pure Appl. Math.}, 52(9):1113--1137, 1999.

\bibitem{ChenLi13}
J.~Chen and Y.~Li.
\newblock Homotopy classes of harmonic maps of the stratified 2-spheres and
  applications to geometric flows.
\newblock {\em Adv. Math.}, 263:357--388, 2014.

\bibitem{FGOT12}
J.~Fan, H.~Gao, T.~Ogawa, and F.~Takahashi.
\newblock A regularity criterion to the biharmonic map heat flow in {$\germ
  R^4$}.
\newblock {\em Math. Nachr.}, 285(16):1963--1968, 2012.

\bibitem{GalPoh02}
V.A. Galaktionov and S.I. Pohozaev.
\newblock Existence and blow-up for higher-order semilinear parabolic
  equations: majorizing order-preserving operators.
\newblock {\em Indiana Univ. Math. J.}, 51(6):1321--1338, 2002.

\bibitem{Gastel02}
A.~Gastel.
\newblock Singularities of first kind in the harmonic map and {Y}ang-{M}ills
  heat flows.
\newblock {\em Math. Z.}, 242(1):47--62, 2002.

\bibitem{Gastel06}
A.~Gastel.
\newblock The extrinsic polyharmonic map heat flow in the critical dimension.
\newblock {\em Adv. Geom.}, 6(4):501--521, 2006.

\bibitem{GastelScheven09}
A.~Gastel and C.~Scheven.
\newblock Regularity of polyharmonic maps in the critical dimension.
\newblock {\em Comm. Anal. Geom.}, 17(2):185--226, 2009.

\bibitem{GastelZorn12}
A.~Gastel and F.~Zorn.
\newblock Biharmonic maps of cohomogeneity one between spheres.
\newblock {\em J. Math. Anal. Appl.}, 387(1):384--399, 2012.

\bibitem{GSZG09}
P.~Goldstein, P.~Strzelecki, and A.~Zatorska-Goldstein.
\newblock On polyharmonic maps into spheres in the critical dimension.
\newblock {\em Ann. Inst. H. Poincar\'e Anal. Non Lin\'eaire},
  26(4):1387--1405, 2009.

\bibitem{Grotowski91}
J.F. Grotowski.
\newblock Harmonic map heat flow for axially symmetric data.
\newblock {\em Manuscripta Math.}, 73(2):207--228, 1991.

\bibitem{HHW2014}
J.~Hineman, T.~Huang, and C.-Y. Wang.
\newblock Regularity and uniqueness of a class of biharmonic map heat flows.
\newblock {\em Calc. Var. Partial Differential Equations}, 50(3-4):491--524,
  2014.

\bibitem{KarcherWood84}
H.~Karcher and J.C. Wood.
\newblock Nonexistence results and growth properties for harmonic maps and
  forms.
\newblock {\em J. Reine Angew. Math.}, 353:165--180, 1984.

\bibitem{Ku08}
Y.~Ku.
\newblock Interior and boundary regularity of intrinsic biharmonic maps to
  spheres.
\newblock {\em Pacific J. Math.}, 234(1):43--67, 2008.

\bibitem{KuwertSchatzle01}
E.~Kuwert and R.~Sch{\"a}tzle.
\newblock The {W}illmore flow with small initial energy.
\newblock {\em J. Differential Geom.}, 57(3):409--441, 2001.

\bibitem{LammSmallEnergy04}
T.~Lamm.
\newblock Heat flow for extrinsic biharmonic maps with small initial energy.
\newblock {\em Ann. Global Anal. Geom.}, 26(4):369--384, 2004.

\bibitem{Lamm05}
T.~Lamm.
\newblock Biharmonic map heat flow into manifolds of nonpositive curvature.
\newblock {\em Calc. Var. Partial Differential Equations}, 22(4):421--445,
  2005.

\bibitem{LammWang09}
T.~Lamm and C.-Y. Wang.
\newblock Boundary regularity for polyharmonic maps in the critical dimension.
\newblock {\em Adv. Calc. Var.}, 2(1):1--16, 2009.

\bibitem{BiharmonicNeckAnalysis}
L.~Liu and H.~Yin.
\newblock {Neck analysis for biharmonic maps}.
\newblock {\em pre-print}, 2013.
\newblock {arXiv:1312.4600 [math.AP]}.

\bibitem{BiharmonicFiniteTimeBlowup}
L.~Liu and H.~Yin.
\newblock {On the finite time blow-up of biharmonic map flow in dimension
  four}.
\newblock {\em pre-print}, 2014.
\newblock {arXiv:1401.6274 [math.AP]}.

\bibitem{MR2013}
S.~Montaldo and A.~Ratto.
\newblock A general approach to equivariant biharmonic maps.
\newblock {\em Mediterr. J. Math.}, 10(2):1127--1139, 2013.

\bibitem{IntrinsicBlowupBehaviour}
R.~Moser.
\newblock The blowup behavior of the biharmonic map heat flow in four
  dimensions.
\newblock {\em IMRP Int. Math. Res. Pap.}, (7):351--402, 2005.

\bibitem{WeakIntrinsic}
R.~Moser.
\newblock Weak solutions of a biharmonic map heat flow.
\newblock {\em Adv. Calc. Var.}, 2(1):73--92, 2009.

\bibitem{Perko91}
L.~Perko.
\newblock {\em Differential equations and dynamical systems}, volume~7 of {\em
  Texts in Applied Mathematics}.
\newblock Springer-Verlag, New York, 1991.

\bibitem{QingTian97}
J.~Qing and G.~Tian.
\newblock Bubbling of the heat flows for harmonic maps from surfaces.
\newblock {\em Comm. Pure Appl. Math.}, 50(4):295--310, 1997.

\bibitem{PierreSchweyer13}
P.~Rapha{\"e}l and R.~Schweyer.
\newblock Stable blowup dynamics for the 1-corotational energy critical
  harmonic heat flow.
\newblock {\em Comm. Pure Appl. Math.}, 66(3):414--480, 2013.

\bibitem{PolyharmonicUniquenessReverseBubbling}
M.~Rupflin.
\newblock Uniqueness for the heat flow for extrinsic polyharmonic maps in the
  critical dimension.
\newblock {\em Comm. Partial Differential Equations}, 36(7):1118--1144, 2011.

\bibitem{Sampson78}
J.H. Sampson.
\newblock Some properties and applications of harmonic mappings.
\newblock {\em Ann. Sci. \'Ecole Norm. Sup. (4)}, 11(2):211--228, 1978.

\bibitem{Strzelecki03}
P.~Strzelecki.
\newblock On biharmonic maps and their generalizations.
\newblock {\em Calc. Var. Partial Differential Equations}, 18(4):401--432,
  2003.

\bibitem{Wang04Arbitrary}
C.-Y. Wang.
\newblock Biharmonic maps from {$\mathbf{R}^4$} into a {R}iemannian manifold.
\newblock {\em Math. Z.}, 247(1):65--87, 2004.

\bibitem{Wang04Spheres}
C.-Y. Wang.
\newblock Remarks on biharmonic maps into spheres.
\newblock {\em Calc. Var. Partial Differential Equations}, 21(3):221--242,
  2004.

\bibitem{Wang04Stationary}
C.-Y. Wang.
\newblock Stationary biharmonic maps from {$\mathbb{R}^m$} into a {R}iemannian
  manifold.
\newblock {\em Comm. Pure Appl. Math.}, 57(4):419--444, 2004.

\bibitem{Chang07}
C.-Y. Wang.
\newblock Heat flow of biharmonic maps in dimensions four and its application.
\newblock {\em Pure Appl. Math. Q.}, 3(2, part 1):595--613, 2007.

\bibitem{BiharmonicRough}
C.-Y. Wang.
\newblock Well-posedness for the heat flow of biharmonic maps with rough
  initial data.
\newblock {\em J. Geom. Anal.}, 22(1):223--243, 2012.

\bibitem{WOY2014}
Z.-P. Wang, Y.-L. Ou, and H.-C. Yang.
\newblock Biharmonic maps from a {$2$}-sphere.
\newblock {\em J. Geom. Phys.}, 77:86--96, 2014.

\bibitem{XuYang02}
X.~Xu and P.C. Yang.
\newblock Conformal energy in four dimension.
\newblock {\em Math. Ann.}, 324(4):731--742, 2002.

\bibitem{ZornThesis}
F.~Zorn.
\newblock {\"{A}quivariante biharmonische Abbildungen}.
\newblock {\em Dissertation, Universit{\"{a}}t Duisburg-Essen}, 2013.

\end{thebibliography}
\end{document}